\documentclass[reqno,fleqn]{amsart}
\usepackage{amsfonts,amsthm,amssymb}
\usepackage{lipsum}
\usepackage[utf8]{inputenc}
\usepackage[english]{babel}

\usepackage{lmodern,textcomp}
\usepackage{amssymb}
\usepackage{amsthm}
\usepackage{amsmath}
\usepackage{mathtools}
\usepackage{latexsym}
\usepackage{tikz}
\usepackage{fancyvrb}
\usepackage{epsfig}
\usepackage{graphicx}
\usepackage{amsmath}
\usepackage{latexsym}
\usepackage{tikz}
\usepackage{rotating}
\usetikzlibrary{positioning}
\usepackage{subcaption}
\usepackage{datetime}
\newtheorem{theorem}{Theorem}[section]

\newtheorem{lemma}[theorem]{Lemma}

\newtheorem{cor}[theorem]{Corollary}
\newtheorem{definition}[theorem]{Definition}

\newtheorem{oprm}[theorem]{Open Problem}
\newtheorem{rem}[theorem]{Remark}

\title[$k$-Distance Magic Labeling and Long Brush  Graphs]{$k$-Distance Magic Labeling and Long Brush Graphs}
\author{\sc V. Vilfred Kamalappan } 
\address{Department of Mathematics, ~Central ~University ~of~ Kerala, Kasaragod, ~India.}

\email{vilfredkamalv@cukerala.ac.in}


\date{}

\begin{document}

\begin{abstract} We define a labeling $f:$ $V(G)$ $\rightarrow$ $\{1, 2, \ldots, n\}$ on a graph $G$ of order $n \geq 3$ as a \emph{$k$-distance magic} ($k$-DM) if $\sum_{w\in \partial N_k(u)}{ f(w)}$ is a constant and independent of $u\in V(G)$ where $\partial N_k(u)$ = $\{v\in V(G): d(u, v) = k\}$, $k\in\mathbb{N}$. Graph $G$ is called a \emph{$k$-DM} if it has a $k$-DM labeling(L). $k$-DML is a generalization of DML or $\Sigma$-labeling of graphs defined by Vilfred. Long Brush, denoted by $LP_{n, m}$, is a graph with vertex set $\{u_1, u_2, . . . , u_n,$ $v_1, v_2, . . . , v_{m}\}$, a path $P_n$ = $u_1$ $u_2$ . . . $u_n$ and edge set $E(P_n)$ $\cup$ $\{u_1v_i:$ $i$ = 1 to $m\}$ $\cup$ $E(<v_1, v_2, . . . , v_{m}>)$ where $<v_1,v_2,...,v_m>$ represents the induced subgraph on $\{v_1, v_2, . . . , v_m\}$, $m+n \geq 3$ and $m,n\in\mathbb{N}$. Long Brush graphs are used to study existence of $k$-DM graphs. In this paper, using partition techniques, we obtain families of $k$-DM graphs and prove that $(i)$ For $k,n \geq 3$, $m \geq 2$ and $k,m,n\in\mathbb{N}$, $LP_{n,m}$ is $k$-DM if and only if $m(m-1) \leq 2n$ and $k$ = $n$; (ii) For every $k\in\mathbb{N}_0$ and a given $m \geq 2$, $LP_{\frac{m(m-1)}{2}+k, m}$ is a $(\frac{m(m-1)}{2}+k)$-DM graph; (iii) For $m \geq 3$, $LP_{1,m}$ = $K_1(u_1)+(K_{m_1} \cup K_{m_2} \cup  ... \cup K_{m_x})$, $x \geq 2$, $1 \leq m_1 \leq m_2 \leq ... \leq m_x$, $m_1+m_2+...+m_x$ = $m$, $m_1+m_2 \geq 3$  and $m_1,m_2,...,m_x,x\in\mathbb{N}$, $LP_{1,m}$ is 2-DM if and only if $u_1$ is assigned with a suitable $j$ and $J_{m+1}\setminus \{j\}$ is partitioned into $x$ constant sum partites of orders $m_1,m_2,...,m_x$, $1 \leq j \leq m+1$; (iv) For $m \geq 2$ if $LP_{2,m}$ contains two pendant vertices, then $LP_{2,m}$ is not a $2$-DM graph; (v) For $m \geq 2$ and $n \geq 3$, if $LP_{n,m}$ contains three pendant vertices, then $LP_{n,m}$ is not a $2$-DM graph; and (vi) for $m_1$ = 1 to 22, we obtain all possible values of $m$ for which $LP_{1, m}$ = $u_1 + (K_{m_1} \cup K_{m_2})$ is 2-DM, $m_1 \leq m_2$, $m = m_1+m_2 \geq 3$ and $m_1,m_2\in\mathbb{N}$.  
\end{abstract}

\maketitle

\section{Introduction}
A vertex labelling of a graph $G$ is an assignment of labels to the vertices of $G$, satisfying certain conditions. More than 200 types of labellings were defined and studied since 1960, by which labelling graphs has become a multidimensional problem. Despite the large number of papers, there are relatively few general results or methods on graph labellings. Labelled graphs serve as useful models for different applications such as Coding Theory, Radar, Astronomy, Circuit Design, X-ray crystallography and Communication Network Addressing \cite{ga21}. 

Partition of numbers seems to be very simple but plays an important role in Combinatorics, Lie theory, representation theory, mathematical physics, and theory of special functions. Euler, Ramajuan, Rademacher and Paul Erodes revealed the beauty and uses of partitions \cite{sa}. In 1987, Vilfred \cite{v87,v96} defined \emph{$\Sigma$-labeling, $\Sigma$-partition} and \emph{$\Sigma$-labeled graphs}. In 2003 \cite{mrs}, the same was independently defined as \emph{1-distance magic vertex labeling} and in a 2009 article \cite{sf} the term `Distance Magic Labeling' was used for the same. The author got motivation to define sigma labeling by noticing two similar situations - the numbers labeled on the faces of a dice and of magic squares \cite{rb}. The sum of the numbers assigned to each pair of opposite faces of a dice is 7 (See page 1 in \cite{v96}). Corresponding to an $n \times n$ magic square with row sum, say $M$, if we form a complete multipartite graph with each row of the square representing a partite set and if we label each vertex with the corresponding integers in the magic squares, then we find that the sum of labels of all vertices in the neighborhood set for each vertex is the same, equal to $(n-1)M$ (See page 97 in \cite{v96}). Since construction of Magic Squares motivated to define $\Sigma$-labeling, the author feels that it is good to use the term `Distance Magic Labeling' for $\Sigma$-labeling and for $1$-Distance Magic labeling.
					
\begin{definition} A labeling $f:$ $V(G)$ $\rightarrow$ $\{1, 2, \ldots, n\}$ on a graph $G$ of order $n \geq 3$ is called a \emph{$k$-distance magic labeling ($k$-DML)} if $\sum_{w\in \partial N_k(u)}{ f(w)}$ is a constant and independent of $u\in V(G)$ where $\partial N_k(u)$ = $\{v\in V(G): d(u, v) = k\}$, $k\in\mathbb{N}$. 
					
Graph $G$ is called a \emph{$k$-distance magic ($k$-DM)} if it has a $k$-DML. 
\end{definition} 
		
\begin{definition} \emph{Long Brush} is defined as a graph $G$ with $V(G) = \{u_1, u_2, . . . , u_n,$ $v_1, v_2,  . . . , v_m\}$, a path $P_n$ = $u_1 u_2 . . . u_n$ and edge set $E(P_n)$ $\cup$ $E(<v_1,v_2,...,v_m>)$ $\cup$ $\{u_1v_i:~ i$ = 1 to $m\}$ where $<v_1,v_2,...,v_m>$ represents the induced subgraph on $\{v_1, v_2, . . . , v_m\}$, $m+n \geq 3$ and $m,n\in\mathbb{N}$. We denote this graph by $LP_{n,m}$. 
\end{definition}
  
 In this paper, we define $k$-Distance Magic ($k$-DM) labeling, $k$-DM graphs and Long Brush graphs $LP_{n, m}$. Using partition techniques, we obtain families of $k$-DM graphs. $k$-DM labeling is a generalization of $1$-DM labeling or DM labeling or $\Sigma$-labeling of graphs defined by Vilfred  \cite{v87}-\cite{v05}. This paper contains 4 sections. Section 1 presents basic defintions and results which are required in the subsequent sections. In Section 2 results on $k$-DML of union of paths and union of cycles are presented. In Section 3, we define Long Brush graphs $LP_{n,m}$ and prove that (i) For $k,n \geq 3$, $m \geq 2$ and $k,m,n\in\mathbb{N}$, $LP_{n,m}$ is $k$-DM if and only if $m(m-1) \leq 2n$ and $k$ = $n$; (ii) For a given $m \geq 2$ and for every $k\in\mathbb{N}_0$, $LP_{\frac{m(m-1)}{2}+k, m}$ is a $(\frac{m(m-1)}{2}+k)$-DM graph. In Sections 4, we use constant sum bipartition of $J_n$ = $\{1,2,...,n\}$ to derive results on 2-DML of $LP_{1, m}$ and prove that (i) For $m \geq 3$, $x \geq 2$,  $LP_{1,m}$ = $K_1(u_1)+(K_{m_1} \cup K_{m_2} \cup  ... \cup K_{m_x})$, $1 \leq m_1 \leq m_2 \leq$ ... $\leq m_x$, $m_1+m_2+...+m_x$ = $m$, $m_1+m_2 \geq 3$  and $m_1,m_2,...,m_x,x\in\mathbb{N}$, $LP_{1,m}$ is 2-DM if and only if $u_1$ is assigned with a suitable $j$ and $J_{m+1}\setminus \{j\}$ is partitioned into $x$ constant sum partites of orders $m_1,m_2,...,m_x$, $1 \leq j \leq m+1$; (ii) For $m \geq 2$ if $LP_{2,m}$ contains two pendant vertices, then $LP_{2,m}$ is not a $2$-DM graph; (iii) For $m \geq 2$ and $n \geq 3$, if $LP_{n,m}$ contains three pendant vertices, then $LP_{n,m}$ is not a $2$-DM graph; and (iv) for $m_1$ = 1 to 22, we calculate all possible values of $m$ for which $LP_{1, m}$ = $u_1 + (K_{m_1} \cup K_{m_2})$ is 2-DM, $m_1 \leq m_1$, $m = m_1+m_2 \geq 3$ and $m_1,m_2\in\mathbb{N}$.		
 
 Search for a more general case of DML is the motivation to define $k$-DML, $k\in\mathbb{N}$. Through out this paper, we consider finite undirected simple graphs and for all basic ideas in graph theory, we follow \cite{dw01}.				

\begin{definition}
	Vertices $u$ and $v$ of a connected graph $G$ are said to be {\em anti-podal} if their distance $d(u, v)$ = $dia(G)$, the diameter of $G$.
\end{definition}

A \emph{necessary condition for a graph $G$ to be a $k$-distance magic} is that $G$ contains at least one component of diameter at least $k$, and at least two distinct $k$-distance neighbourhoods, $k\in\mathbb{N}$. These conditions need not be a sufficient one.

The following result is useful to identify certain graphs as non-$k$-distance magic (non-$k$-DM) graphs.

\begin{lemma} \quad \label{a1} {\rm For $k \in \mathbb{N}$ and  graph $G$, if $u,v\in V(G)$ such that 
		$$ | \partial N_k(u) \setminus \partial N_k(v)| = 1 = |\partial N_k(v) \setminus \partial N_k(u)|,$$ 
then $G$ is not a $k$-DM graph.}
\end{lemma}

\begin{proof}
If possible, let $G$ be $k$-DM and $f$ be a $k$-DML of $G$. Let $\partial N_k(u) \setminus \partial N_k(v) = \lbrace x \rbrace$ and $ \partial N_k(v) \setminus \partial N_k(u)$ = $\lbrace y \rbrace$, $x,y\in V(G)$. This implies,
$\partial N_k(u)$ = $\partial N_k(v) \cup \lbrace x \rbrace$, $\partial N_k(v)$ = $\partial N_k(u) \cup \lbrace y \rbrace$ and $x \neq y$.

Since $f$ is a $k$-DML of $G$, we get, 
$$\sum_{w\in \partial N_k(u)}{ f(w)} = \sum_{w\in \partial N_k(v)}{ f(w)}.$$
$$\Rightarrow f(x) + \sum_{w\in \partial N_k(v)}{ f(w)} = f(y) + \sum_{w\in \partial N_k(u)}{ f(w)}.$$ 
$\Rightarrow f(x) = f(y)$, which is a contradiction since $f$ is a $k$-DML of $G$, $x \neq y$ and $x,y \in V(G).$ Hence, we get the result.
\end{proof}

\section{$k$-DML of union of paths and union of cycles}

In this section, we study $k$-DML of union of paths and union of cycles, $k\in\mathbb{N}$.

\begin{theorem} \quad \label{a2} {\rm  For $n \geq 2$ and $k,n\in\mathbb{N}$, path $P_n$ is $k$-distance magic if and only if $k = 1$ and $n = 3$.}
\end{theorem}
\begin{proof} Let $P_n = u_1 u_2 \ldots u_n$, $n \geq 2$. Diameter of $P_n$ is $n-1$ and so for $k$-DML of $P_n$, $k \leq n-1$, $k\in\mathbb{N}$. Let $n > k$ and $k,n\in\mathbb{N}$. 

If possible, let $P_n$ be a $k$-DM graph and $f$ be a $k$-DML of $P_n$. 
$$\therefore ~ \sum_{u\in \partial N_k(u_1)}{ f(u)} = \sum_{u\in \partial N_k(u_{n})}{ f(u)}.~ \Rightarrow ~ f(u_{k+1}) = f(u_{n-k}).$$
This is possible only when $f(u_{k+1}) = f(u_{n-k})$, $k\in\mathbb{N}$. That is when $n$ = $2k+1$, $k\in\mathbb{N}$. When $k$ = 1, $n$ = 3 and in this case, $P_n$ = $P_3$ = $u_1 u_2 u_3$ is DM and its DML is given by $\{f(u_1), f(u_3)\}$ = $\{1, 2\}$ and $f(u_2)$ = 3. Graph $P_3$ and its DML are given in Figure 1. See Figure 1.

	%Fig.1,2
	%%%%%%%%%%%%%%%%%%%%%%%%%%%
	\begin{center} 
	\begin{tikzpicture}  \label{f2}
	\node (u1) at (10,1.5)  [circle,draw,scale=.7]{$u_1$};
	\node (u2) at (11,0)  [circle,draw,scale=.7]{$u_2$};
	\node (u3) at (12,1.5)  [circle,draw,scale=.7] {$u_3$};
	
	\draw (u2) -- (u1);
	\draw (u2) -- (u3);
	
	%Fig 2	
	\node (y1) at (15,1.5)  [circle,draw,scale=.7]{1};
	\node (y2) at (17,1.5)  [circle,draw,scale=.7] {2};
	\node (x1) at (16,0)  [circle,draw,scale=.7]{3};
	
	\draw (x1) -- (y1);
	\draw (x1) -- (y2);
	
	\end{tikzpicture}

\vspace{.1cm}		
 \small{$LP_{1, 2} = P_3$ \hspace{3cm}  $P_3$ with DML} 

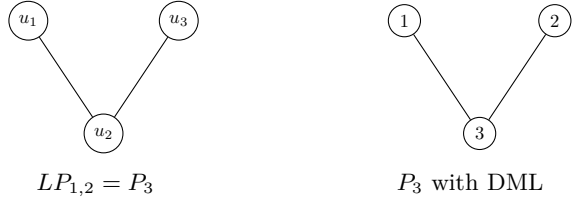
 \small{ Figure 1.} 
\end{center}
%%%%%%%%%%%%%%%%%%%%%%%%%%%%%%

When $k > 1$ and $n$ = $2k+1$, $n \geq 5$, $P_n$ = $P_{2k+1}$ = $u_1 u_2 . . . u_{2k+1}$, $\partial N_k(u_{k+2})$ = $\{u_2\}$, $\partial N_k(u_{k+3})$ = $\{u_3\}$ and 

$$\sum_{u\in \partial N_k(u_{k+2})}{ f(u)} = \sum_{u\in \partial N_k(u_{k+3})}{ f(u)}.~ \Rightarrow ~ f(u_{2}) = f(u_{3})$$
which is a contradiction and thereby $k$-DML doesn't exist to $P_n$ when $k > 1$ and $n$ = $2k+1$. Thus from the above cases, we get the result. 
\end{proof}

\begin{theorem} \quad \label{a3} {\rm For $n_i \geq 2$, $1 \leq i \leq m$ and $k,m,n_i\in\mathbb{N}$, graph $\cup_{i=1}^m P_{n_i}$ is $k$-DM if and only if $n_i$ = 3 and $m$ = 1 = $k$.}
\end{theorem}
\begin{proof} If $\cup_{i=1}^m P_{n_i}$ is $k$-DM, then by applying the same proof technique of Theorem \ref{a2} on each component $P_{n_i}$, we get $n_i$ = 3 and $k$ = 1 for each $i$, $1 \leq i \leq m$. Let $P_{n_i} = u^{(i)}_{1} u^{(i)}_{2} u^{(i)}_{3}$, $n_i = 3$, $k$ = 1. Also, $n$, the biggest number among the vertex labels, is assigned to a vertex, say $u^{(i)}_{j}$ such that $\partial N_k(u^{(i)}_{l})$ = $\partial N_1(u^{(i)}_{l})$ = $N(u^{(i)}_{l})$ = $\{u^{(i)}_{j}\}$ for some $u^{(i)}_{l}$ and $f(u^{(i)}_{j})$ = $n$, $k$ = 1, $1 \leq i \leq m$ and $1 \leq j,l \leq 3$. This implies, $j$ = 2, $l$ = 1 or 3 and so $f(u^{(i)}_{2})$ = $n$, $\forall i$, $1 \leq i \leq m$. This is possible only when $m$ = 1 since $f$ is a DML. Hence we get the result.
\end{proof}

\begin{theorem} {\rm  {[Proposition 2.2.8 in \cite{v96}]}  \quad \label{a4} Any component of a DM graph which is a cycle must be of length 4. \hfill $\Box$} 
\end{theorem}
\begin{cor} {\rm \quad \label{a5} For $n_i \geq 3$, $1 \leq i \leq m$ and $m,n_i\in\mathbb{N}$, graph $\cup_{i=1}^m C_{n_i}$ is DM if and only if $n_i$ = 4 for every $i$. } 
\end{cor}
\begin{proof}\quad Here, we present a proof different from that of Proposition 2.2.8 in \cite{v96}. 

 Let $C_{n_j} = (u^{(j)}_{1}, u^{(j)}_{2}, . . . , u^{(j)}_{n_j})$, $1 \leq j \leq m$. Let $G$ = $\cup_{i=1}^m C_{n_i}$ be DM and $f$ be a DML of $G$. Our aim is to prove that $n_i$ = 4 for all $i$, $1 \leq i \leq m$. 
 
 For $1 \leq j \leq m$, in $C_{n_j}$, 
$$N(u^{(j)}_{2}) = \{u^{(j)}_{1}, u^{(j)}_{3}\} ~and~ N(u^{(j)}_{n_j}) = \{u^{(j)}_{1}, u^{(j)}_{n_j-1}\}.$$
$$\Rightarrow~ f(u^{(j)}_{1}) + f(u^{(j)}_{3}) = f(u^{(j)}_{1}) + f(u^{(j)}_{n_j-1}).~ \Rightarrow~ f(u^{(j)}_{3}) = f(u^{(j)}_{n_j-1}).$$
$\Rightarrow$ $n_j-1$ = 3 since $f$ is a DML on $G$. This implies, $n_j$ = 4. 

This is true for every $C_{n_j}$, $1 \leq j \leq m$. Thus $n_j$ = 4 for every $C_{n_j}$, $1 \leq j \leq m$.

Conversely, let $G$ = $m. C_4$, $m\in\mathbb{N}$. Our aim is to prove that $G$ is DM.

$G$ is union of $m$ disjoint copies of $C_4$. Let $C^{(j)}_4$ = $(u^{(j)}_{1}, u^{(j)}_{2}, u^{(j)}_{3}, u^{(j)}_{4})$ be the $j^{th}$ copy of $C_4$ in $G$, $1 \leq j \leq m$. Let  
$$f: \{\{u^{(j)}_{1}, u^{(j)}_{3}\}, \{u^{(j)}_{2}, u^{(j)}_{4}\}: j = 1~to~m\} \to \{\{1, 4m\}, \{2, 4m-1\}, ..., \{2m, 2m+1\}\}$$ be a bijective mapping. Clearly, 
$$\sum_{v\in N(u^{(j)}_{1})}f(v) = \sum_{v\in N(u^{(j)}_{4})}f(v) = 4m+1~ since$$
$$f(u^{(j)}_{2}) + f(u^{(j)}_{4}) = k+(4m-k+1) ~and ~f(u^{(j)}_{1}) + f(u^{(j)}_{3}) = l+(4m-l+1)$$ for some $l$ and $k$, $1 \leq k,l \leq 2m$ and $j$ = 1 to $m$.

$\Rightarrow$ $f$ is a DML on $G$ and hence $G$ is a DM graph. 

Hence we get the result.
\end{proof}

\begin{cor} {\rm \quad \label{a6} The number of distinct DMLs of the labeled graph $G$ = $m. C_{4}$ is $2^{2m}(2m)!$, $m\in\mathbb{N}$. } 
\end{cor}
\begin{proof}\quad The number of distinct DMLs of the labeled graph $m. C_{4}$ is same as the number of distinct bijective mappings $f$ as defined in the proof of Corollary \ref{a5} since constant sum partition of $\{1, 2, . . . , 4m\}$ with each partite of order 2 is $\{1, 4m\}$, $\{2, 4m-1\}$, . . . , $\{2m, 2m+1\}$ only. 

Here, after selecting a pair of labels, say $\{i, 4m-i+1\}$ to a pair of non-adjacent vertices, say, $\{u^{(j)}_{1}, u^{(j)}_{3}\}$ in $C_{n_j}$ = $C^{(j)}_{4}$ = $(u^{(j)}_{1}, u^{(j)}_{2}, u^{(j)}_{3}, u^{(j)}_{4})$, there are two possible ways to label vertices $u^{(j)}_{1}$ and $u^{(j)}_{3}$, $1 \leq j \leq m$. And number of bijective mappings from $\{a_1, a_2, . . . , a_{2m}\}$ $\to$ $\{b_1, b_2, . . . , b_{2m}\}$ is $(2m)!$. Hence the total number of DMLs on the labeled graph $G$ = $m.C_4$ is $2^{2m}(2m)!$, $m\in\mathbb{N}$.
\end{proof}

\begin{theorem}  \quad \label{a7} {\rm Let $k \geq 2$, $n \geq 3$ and $k,n\in\mathbb{N}$. Then $C_n$ is $k$-DM if and only if $n = 4k$.}
\end{theorem}
\begin{proof}\quad Let $C_n = (u_1, u_2, . . . , u_{n})$, $n \geq 3$. 

When $n = 4k$, $k \geq 2$ and $k,n\in\mathbb{N}$, $C_n = C_{4k} = (u_1, u_2, . . . , u_{4k})$. Let 
$$f: \{\{u_{i}, u_{2k+i}\}: i = 1~to~2k\} \to \{\{j, 4k+1-j\}: j = 1 ~to ~2k\}$$ 
be a bijective mapping. 

Clearly, for $i$ = 1 to $2k$, under subscript modulo $4k$ and $u_0$ = $u_{4k}$,
$$\sum_{v\in \partial N_k(u_{i})}f(v) = f(u_{k+i}) + f(u_{3k+i}) = j + (4k+1-j) ~for ~some~ j,~ 1 \leq j \leq 2k.$$
$$ = 4k+1 = \sum_{v\in \partial N_k(u_{2k+i})}f(v).$$

$\Rightarrow$ $f$ is a $k$-DML on $C_{4k}$ = $C_n$. 

$\Rightarrow$ $C_{4k}$ is a $k$-DM graph, $k \geq 2$ and $k\in\mathbb{N}$.

Conversely, let $C_n$ be a $k$-DM graph and $f$ be a $k$-DML of $C_n$. Our aim is to prove that $n$ = $4k$, $k\in\mathbb{N}$. We consider the following 3 cases of $n$.

\vspace{.1cm}
\noindent
\emph{Case 1.}\quad $n < 2k$, $k \geq 2$ and $k,n\in\mathbb{N}$.

In this case, $k$-DML does't exist since $dia(C_n) < k$ when $n < 2k$.

\vspace{.1cm}
\noindent
\emph{Case 2.}\quad $n = 2k$, $k \geq 2$ and $k\in\mathbb{N}$.

In this case, $C_n = C_{2k} = (u_1, u_2, . . . , u_{2k})$, $dia(C_n)$ = $dia(C_{2k})$ = $k$, $\partial N_k(u_{1})$ = $\{u_{k+1}\}$ and $\partial N_k(u_{k+1})$ = $\{u_{1}\}$ which implies, $f(u_{k+1})$ = $f(u_1)$ since $C_n$ is $k$-DM graph and $f$ is a $k$-DML of $C_n$. This implies, $u_{k+1}$ = $u_1$ which is a contradiction since $u_{k+1}$ $\neq$ $u_1$. Thus, when $n$ = $2k$, $C_n$ is not a $k$-DM graph.

\vspace{.1cm}
\noindent
\emph{Case 3.}\quad $n > 2k$, $k \geq 2$ and $k\in\mathbb{N}$.

In this case, $dia(C_n)$ $\geq$ $k$, $\partial N_k(u_{k+1})$ = $\{u_1, u_{2k+1}\}$ and 

$\partial N_k(u_{n-k+1})$ = $\{u_1, u_{n-2k+1}\}$.

$\Rightarrow$ $f(u_1)$ + $f(u_{n-2k+1})$ = $f(u_1)$ + $f(u_{2k+1})$ since $f$ is a $k$-DML of $C_n$.

$\Rightarrow$ $f(u_{2k+1})$ = $f(u_{n-2k+1})$, $1 \leq n-2k+1 \leq n$.

$\Rightarrow$ $u_{2k+1}$ = $u_{n-2k+1}$ since $f$ is a $k$-DML of $C_n$.

$\Rightarrow$ $2k+1$ = $n-2k+1$, $1 \leq n-2k+1 \leq n$.

$\Rightarrow$ $n$ = $4k$ is the only possibility when $n > 2k$, $k \geq 2$ and $k,n\in\mathbb{N}$. 

Hence we get the result.
\end{proof}

\begin{theorem} {\rm \quad \label{a8} Let $G$ = $\cup_{i=1}^m C_{n_i}$, $k \geq 2$, $n_i \geq 2k$, $\forall i$, $1 \leq i \leq m$ and $m,n_i\in\mathbb{N}$. Then $G$ is $k$-DM if and only if $n_i$ = $4k$ for every $i$, $1 \leq i \leq m$. } 
\end{theorem}
\begin{proof}\quad Let $C_{n_j} = (u^{(j)}_{1}, u^{(j)}_{2}, . . . , u^{(j)}_{n_j})$, $1 \leq j \leq m$. Let $n_j = 4k$ for all $j$, $1 \leq j \leq m$ and $k \geq 2$. Thus $G$ = $\cup_{i=1}^m C_{4k}$. Let $C^{(j)}_{4k}$ = $(u^{(j)}_{1}, u^{(j)}_{2}, . . . , u^{(j)}_{4k})$ be the $j^{th}$ copy of $C_{4k}$ in $G$, $1 \leq j \leq m$ and $j,m\in\mathbb{N}$. Let 
$$f: {\{\{u^{(j)}_{i}, u^{(j)}_{2k+i}\}: i = 1~to~2k~\&~j = 1~to~m\}} \to {\{\{l, 4km+1-l\}: l = 1~to~2km\}}$$  be a bijective mapping. 

Clearly, for $i$ = 1 to $2k$, $1 \leq j \leq m$ and $u^{(j)}_{0}$ = $u^{(j)}_{4k}$, under subscript modulo $4k$, 
$$\sum_{v\in \partial N_k(u^{(j)}_{i})}f(v) = f(u^{(j)}_{k+i}) + f(u^{(j)}_{3k+i}) = l + (4km+1-l) ~for ~some~ l,~ 1 \leq l \leq 2km$$
$$ = 4km+1 = \sum_{v\in \partial N_k(u^{(j)}_{2k+i})}f(v).$$
$\Rightarrow$ $f$ is a $k$-DML on $G$ = $\cup_{i=1}^m C_{n_i}$, $k \geq 2$, $n_i = 4k$, $\forall i$, $1 \leq i \leq m$ and $k,m\in\mathbb{N}$. 

$\Rightarrow$ For $k \geq 2$, $n_i = 4k$, $\forall i$, $1 \leq i \leq m$ and $k,m\in\mathbb{N}$, graph $G = \cup_{i=1}^m C_{n_i}$ is a $k$-DM graph. 

For the converse part, apply the same proof technique of Theorem \ref{a7} on each $C_{n_j}$, we get $n_j$ = $4k$ for all $j$, $1 \leq j \leq m$. Figure 2 shows relative positions of $k$ distance points $u^{(j)}_{i}$, $u^{(j)}_{k+i}$, $u^{(j)}_{2k+i}$, $u^{(j)}_{3k+i}$ of $C_{n_j}$, $1 \leq i \leq k$. See Figure 2. 

Hence the result is proved.
\end{proof}

	%Fig.3	
	\begin{center} 
		\begin{tikzpicture} 
	\node (u1) at (1,1)  [circle,draw,scale=2.6] {};
	\node (u2) at (2.5,-1)  [circle,draw,scale=2.6] {};
	\node (u3) at (-0.5,-1)  [circle,draw,scale=2.6] {};
	\node (u4) at (-2,.5)  [circle,draw,scale=2.6] {};
	
	\node (v1) at (1,1)  {$u^{(j)}_{i}$};
	\node (v2) at (2.5,-1)  {$u^{(j)}_{k+i}$};
	\node (v3) at (-0.5,-1)  {$u^{(j)}_{2k+i}$};
	\node (v4) at (-2,.5)   {$u^{(j)}_{3k+i}$};
	
  \node (v5) at (0,0)  {$1 \leq i \leq k.$};
 	
	\draw (u1) [dashed] -- (u2);
	\draw (u1) [dashed] -- (u4);
	\draw (u2) [dashed] -- (u3);
	\draw (u3) [dashed] -- (u4);
	
	\node (v) at (0,-2)  {Figure 2};
\end{tikzpicture}
\end{center} 

\begin{cor} {\rm \quad \label{a9} The number of distinct $k$-DMLs on the labeled graph $\cup_{i=1}^m C_{4k}$ is $2^{2km}(2km)!$,  $k,m\in\mathbb{N}$. } 
\end{cor}
\begin{proof} \quad Similar to the proof of Corollary \ref{a6}. 
\end{proof}
For $k$ = 1, the above result becomes Corollary \ref{a6}.
 
	\section{Long Brush graphs $LP_{n,m}$ and their $k$-DMLs}
	
	In this section, we study $k$-DML of Long Brush graphs $LP_{n,m}$. Graph $LP_{n,m}$ has vertex set $\{u_1,u_2,...,$ $u_n, v_1,v_2,...,v_m\}$, edge set $E(P_n)$ $\cup$ $E(<v_1,v_2,...,v_m>)$ $\cup$ $\{u_1v_i:~ i$ = 1 to $m\}$ and contains the path $P_n$ = $u_1 u_2 . . . u_n$. In $LP_{n,m}$, the induced subgraph $<u_1, v_1, v_2, . . . , v_m>$ is called the \emph{brush} and the path $P_n$ = $u_1 u_2 . . . u_n$ as the \emph{handle} of the Long Brush $LP_{n,m}$, $m+n \geq 3$ and $m,n\in\mathbb{N}$.

 For $n \geq 2$, in $LP_{n,m}$, $p = |V(LP_{n,m})|$ = $m+n$, $q = |E(LP_{n,m})|$ = $n-1+m$ + $|E(<v_1, v_2,$ $. . . , v_m>)|$, $m+n-1 \leq q \leq m+n-1$ + $\left( {\begin{array}{c}
   $m$ \\
   2  
 \end{array} } \right)$, $dia(LP_{n,m})$ = $n$ and $u_n$ and $v_i$ are pair of antipodal vertices for every $i$, $1 \leq i \leq m$. Moreover, here $q$ = $m+n-1$ when $<v_1, v_2,  . . . , v_m>$ = $\overline{K}_m$ and $q$ = $m+n-1$ + $\left( {\begin{array}{c}
   $m$ \\
   2  
 \end{array} } \right)$ when $<v_1, v_2,  . . . , v_m>$ = $K_m$. 

When $m$ = 1, $LP_{n,m}$ = $LP_{n,1}$ = $P_{n+1}$ = $v_1 u_1 u_2 . . . u_n$ which is DM only when $n+1$ = 3. And for $k \geq 2$, $P_{n+1}$ is not a $k$-DM graph by Theorem \ref{a2}. Thus hereafter while discussing $k$-DML of $LP_{n,m}$, we consider $k,m \geq 2$ and $k,m,n\in\mathbb{N}$. In this study, we consider long brush graphs $LP_{n,m}$ with $k,n \geq 3$ and $m \geq 2$ at first and then we consider the case of $k$ = 2 and $m,n \geq 2$, $k,m,n\in\mathbb{N}$.

\begin{theorem}  \quad \label{b1} {\rm For $k,n \geq 3$ and $m \geq 2$, $LP_{n,m}$ is $k$-DM if and only if $m(m-1) \leq 2n$ and $k$ = $n$, $k,m,n\in\mathbb{N}$.}
\end{theorem}
\begin{proof}\quad  Clearly, for $m \geq 2$ and $k,n \geq 3$, $dia(LP_{n,m})$ = $n$ and hence $k \leq n$ when $LP_{n,m}$ is a $k$-DM graph. 
	
Let $LP_{n,m}$ be $k$-DM and $f$ be a $k$-DML of $LP_{n,m}$, $k \geq 3$.

When $k \geq 3$, in $LP_{n,m}$, $d(v_i, v_j) \leq 2$, $\partial N_k(v_i)$ = $\{u_k\}$ and $\partial N_k(u_{k}) = \{v_1,v_2,...,$ $v_m\}$ if $n < 2k$ and = $\{v_1,v_2,...,v_m\}$ $\cup$ $\{u_{2k}\}$ if $n \geq 2k$, $\forall i, j$, $1 \leq i \leq j \leq m$.  
 $$\Rightarrow~ \sum_{u\in \partial N_k(v_{i})}f(u) = f(u_{k}) =  \sum_{u\in \partial N_k(u_{k})}f(u)~ and ~ \hspace{3.7cm} $$
 \[ \sum_{u\in \partial N_k(u_{k})}f(u) = \left\{ \begin{array}{ll}
		\sum^m_{i=1}{f(v_i)} & \mbox{if $n < 2k$ and $k \geq 3$} \\
			\sum^m_{i=1}{f(v_i)} + f(u_{2k}) & \mbox{if $n \geq 2k$ and $k \geq 3$.}
		\end{array}  \right.  \]
And thereby $u_k$ has to take possible bigger value and $v_1, v_2,  . . . , v_m$ smaller values from $J_{m+n}$ = $\{1, 2, . . . , m+n\}$.
This implies,  
$$\sum^m_{i=1}{f(v_i)} \geq 1+2+...+m = \frac{m(m+1)}{2} ~and ~ \sum^m_{i=1}{f(v_i)} \leq f(u_k),~ k \geq 3.$$

$\Rightarrow$ $\frac{m(m+1)}{2} \leq m+n$ since the possible  maximum value of $f(u_k)$ in $LP_{n,m}$ is $m+n$ = $|V(LP_{n,m})|$. This implies, $m(m-1) \leq 2n$.

Thus, for $k,n \geq 3$ and $m \geq 2$, if $LP_{n,m}$ is $k$-DM, then $k \leq n$ and $m(m-1) \leq 2n$.

To complete the proof, we have to prove that when $k < n$ and $k,n \geq 3$, $LP_{n,m}$ is not $k$-DM and $LP_{n,m}$ is $k$-DM when $k = n$. This is done by the following Claims 1 and 2.

\vspace{.1cm}
\noindent
{\it Claim 1.}\quad $LP_{n,m}$ is not $k$-DM when $k < n$, $k,n \geq 3$ and $m \geq 2$. 

If possible, let $LP_{n,m}$ be $k$-DM and $f$ be a $k$-DML of $LP_{n,m}$, $k < n$, $k,n \geq 3$ and $m \geq 2$. Then using the definition of $k$-DML, we get, 
$$\sum_{v\in \partial N_k(v_{i})}{f(v)} = \sum_{v\in \partial N_k(u_{1})}{f(v)},~ 1 \leq i \leq m.$$
$$\Rightarrow f(u_k) = f(u_{k+1})~ since ~ \partial N_k(v_{i}) = \{u_k\}~ and ~ \partial N_k(u_{1}) = \{u_{k+1}\},$$ $k+1 \leq n$, $k,n \geq 3$, $m \geq 2$ and $1 \leq i \leq m.$ This is a contradiction to our assumption that $f$ is a $k$-DML of $LP_{n,m}$. Hence the claim is true.

\vspace{.1cm}
\noindent
\emph{Claim 2.}\quad For $k,n \geq 3$, $m \geq 2$ and $m(m-1) \leq 2n$, $LP_{n,m}$ is $n$-DM and $k$ = $n$.

It is enough to prove that $LP_{n,m}$ is $n$-DM when $m(m-1) \leq 2n$, $n \geq 3$ and $m \geq 2$. Clearly, $dia(LP_{n,m}) = n$ and hence $\partial N_n(v_{i})$ = $\{u_n\}$ and $\partial N_n(u_{n})$ = $\{v_1, v_2, . . . , v_m\}$. 

Given that  $m(m-1) \leq 2n$ which implies, $\frac{m(m+1)}{2} \leq n+m$, $\frac{m(m+1)}{2}\in\mathbb{N}$. \hfill (a)

For $m \geq 2$, there are two possiblities now, (i) $m+1 = \frac{m(m+1)}{2}$ or (ii) $m+1 < \frac{m(m+1)}{2}$.  Crrespondingly, we consider the following two cases. 
\begin{enumerate}
\item [\rm {\bf Case (i)}]  $m \geq 2$ and $m+1 = \frac{m(m+1)}{2}$, $m\in\mathbb{N}$. 

In this case, we get $m$ = 2 and so $m+n$ = $n+2$  = $|V(LP_{n,2})|$. Now, consider a bijective mapping 
$$f: V(LP_{n,2}) \to \{1, 2, . . . , n+2\} \ni$$  
$$~\{f(v_1), f(v_2)\} = \{n+2-k, k\},~ 1 \leq k \leq n+1$$ 
$$f(u_n) = n+2~ and$$
$$f(\{u_1, u_2, . . . , u_{n-1}\}) = [1, n+1] \setminus \{f(v_1), f(v_2)\}.$$
Clearly, $f$ is an $n$-DML of $LP_{n,2}$ with constant sum $M$ = $n+2$ for $m = 2$ and $n \geq 3$. 
\item [\rm {\bf Case (ii)}] $m \geq 2$ and $m+1 < \frac{m(m+1)}{2}$, $m\in\mathbb{N}$

In this case, $m > 2$ and $m+1 < \frac{m(m+1)}{2} \leq n+m$ using (a). \hfill (b)

Consider a bijective mapping 
$$g: V(LP_{n,m}) \to \{1, 2, . . . , m+n\} \ni$$  
$$~ ~g(v_i) = i ~for ~ i = 1~ to ~m,$$ 
$$g(u_n) = \frac{m(m+1)}{2}~ and$$
$$g(\{u_1, u_2, . . . , u_{n-1}\}) = \{m+1, m+2, . . . , \frac{m(m+1)}{2}-1,$$
$$\hspace{6cm} \frac{m(m+1)}{2}+1, \frac{m(m+1)}{2}+2, . . . , m+n\}.$$

Clearly, $g$ is an $n$-DML of $LP_{n,m}$ with constant sum $\frac{m(m+1)}{2}$ for $m > 2$ and $n \geq 3$. Thereby the claim is true in this case.
\end{enumerate}
 Thus for $m \geq 2$, $k,n \geq 3$ and $k,m,n\in\mathbb{N}$, when $m(m-1) \leq n+m$, $LP_{n,m}$ is $k$-DM and $k$ = $n$.

Hence, we get, for $m \geq 2$, $k,n \geq 3$ and $k,m,n\in\mathbb{N}$, $LP_{n,m}$ is $k$-DM if and only if $k$ = $n$ and $m(m-1) \leq n+m$. 
\end{proof}

Structure of graph $LP_{n,m}$ depends mainly on the structure of $<v_1,v_2,...,v_m>$. The follwing theorems present results on 2-DM of $LP_{n,m}$ when $<v_1,v_2,...,v_m>$ = $K_m$ as well as it contains either one isolated vertex when $n$ = 2  or two isolated vertices when $n \geq 3$.

\begin{theorem}  \quad \label{b2} {\rm Let $m,n \geq 2$ and $LP_{n,m}$ be 2-DM. Then  

(i) if $<v_1, v_2, . . . , v_m>$ = $K_m$, then $m(m-1) \leq n+m$;  

(ii) if $<v_1, v_2, . . . , v_m>$ $\neq$ $K_m$, then $\forall$ $v_i,v_j \ni d(v_i, v_j)$ = 2 and $1 \leq i < j \leq m$, 
$$\sum_{v_x\in N[v_i]}{f(v_x)} = \sum_{v_x\in N[v_j]}{f(v_x)}.$$}
\end{theorem}
\begin{proof}\quad  
\begin{enumerate}
	\item  [\rm (i)] When  $k$ = 2, if $<v_1, v_2, . . . , v_m>$ = $K_m$, then similar to the Case 1 in the proof of Theorem \ref{b1}, by applying $k$ = 2, we get,   

$$\hspace{1cm} \sum^m_{i=1}{f(v_i)} \geq 1+2+...+m = \frac{m(m+1)}{2} ~and ~\partial N_2(v_j) = \{u_2\},~ 1 \leq j \leq m.$$
Also, for $1 \leq i \leq m$ and $k = 2$, 
$$\sum_{u\in \partial N_2(v_{i})}{f(u)} = f(u_{2})~ and \hspace{4cm}$$

 \[ f(u_2) = \left\{ \begin{array}{ll}
		\sum^m_{i=1}{f(v_i)} & \mbox{if $n < 4$} \\
			\sum^m_{i=1}{f(v_i)} + f(u_{4}) & \mbox{if $n \geq 4$.}
		\end{array}  \right.  \]
This implies, $f(u_2) \geq 	\sum^m_{i=1}{f(v_i)} \geq \frac{m(m+1)}{2}$. But $f(u_2) \leq m+n$ = $p$ = $|V(LP_{n,m})|$. This implies, $\frac{m(m+1)}{2} \leq$ $m+n$ which implies, $m(m-1) \leq 2n$.

 Thus, for $k$ = 2, $m,n \geq 2$ and $<v_1, v_2, . . . , v_m>$ = $K_m$, if $LP_{n,m}$ is $2$-DM, then $m(m-1) \leq 2n$.

\item [\rm (ii)] When $k$ = 2, $m,n \geq 2$, $<v_1, v_2, . . . , v_m>$ $\neq$ $K_m$ and $LP_{n,m}$ is $2$-DM with $f$ as a $2$-DML, then there exists $v_i,v_j \ni v_iv_j\notin E(<v_1, v_2, . . . , v_m>)$, $1 \leq i < j \leq m$. Thus, for $1 \leq i < j \leq m$ and $v_i,v_j\in V(LP_{n,m})$ $\ni$ $v_iv_j\notin E(<v_1, v_2, . . . , v_m>)$, in $LP_{n,m}$,
$$\partial N_2(v_i) = \{u_2\} \cup \{ v_x: d(v_i,v_x) = 2, ~1 \leq x \leq m\}; \hspace{2cm} (1)$$
$$\partial N_2(v_j) = \{u_2\} \cup \{ v_y: d(v_j,v_y) = 2, ~1 \leq y \leq m\};  \hspace{2cm} (2)$$ 
\[  \hspace{1cm} \partial N_2(u_2) = \left\{ \begin{array}{ll}
		\{v_i: i$ = 1 to $m\}  & \mbox{if $n < 4$} \\
			\{u_{4}, v_i: i$ = 1 to $m\} & \mbox{if $n \geq 4.$\hspace{2.1cm} (3)}
		\end{array}  \right.  \]

\[ \Rightarrow  \sum_{v\in\partial N_2(u_2)}{f(v)} = \left\{ \begin{array}{ll} 
		\sum^m_{i=1}{f(v_i)}  & \mbox{if $n < 4$} \\
			\sum^m_{i=1}{f(v_i)} + f(u_{4}) & \mbox{if $n \geq 4.$\hspace{2cm} (4)}
		\end{array}  \right.  \]
Since ~$f$ is a 2-DML of $LP_{n,m}$, we get,
$$ \sum_{v\in \partial N_2(v_i)}{f(v)} = \sum_{v\in \partial N_2(v_j)}{f(v)} = \sum_{v\in\partial N_2(u_2)}{f(v)},$$
\hfill $1 \leq i < j \leq m$ and $d(v_i, v_j) = 2.$
$$\Rightarrow f(u_2) + \sum_{d(v_i,v_x) = 2,~1 \leq x \leq m}{f(v_x)} = f(u_2) + \sum_{d(v_j,v_x) = 2,~1 \leq x \leq m}{f(v_x)}.$$ 

$$\Rightarrow  \sum_{d(v_i,v_x) = 2,~1 \leq x \leq m}{f(v_x)} =  \sum_{d(v_j,v_x) = 2,~1 \leq x \leq m}{f(v_x)},$$
\hfill $1 \leq i < j \leq m$ and $d(v_i, v_j) = 2.$ 
$$\Rightarrow  \sum^m_{x=1}{f(v_x)} - \sum_{v_x\in N[v_i]}{f(v_x)} =  \sum^m_{x=1}{f(v_x)} - \sum_{v_x\in N[v_j]}{f(v_x)},$$
\hfill $1 \leq i < j \leq m$ and $d(v_i, v_j) = 2.$ 
$$\Rightarrow  \sum_{v_x\in N[v_i]}{f(v_x)} =  \sum_{v_x\in N[v_j]}{f(v_x)}, ~ 1 \leq i < j \leq m~ and~ d(v_i, v_j) = 2.$$
\end{enumerate}	
Hence we get the result.
\end{proof}

\begin{theorem}  \quad \label{b23} {\rm 
\begin{enumerate} \item [\rm (i)]\quad For $m \geq 2$, in $LP_{2,m}$, if $<v_1, v_2, . . . , v_m>$ = $K_1(v_1)$ $\cup$ $<v_2, v_3, . . . , v_m>$, then $LP_{2,m}$ is not a $2$-DM graph.
		
i.e., $LP_{2,m}$ is not a $2$-DM graph when it contains one pendant vertex.
\item [\rm (ii)]  For $m \geq 2$ and $n \geq 3$, in $LP_{n,m}$, if $<v_1, v_2, . . . , v_m>$ = $K_1(v_1) \cup K_1(v_2)$ $\cup$ $<v_3, . . . , v_m>$, then $LP_{n,m}$ is not a $2$-DM graph.

i.e., $LP_{n,m}$  is not a $2$-DM graph when it contains two pendant vertices.
\end{enumerate}}
\end{theorem}
\begin{proof}\quad \begin{enumerate} \item [\rm (i)]\quad Here, $v_1$ is a pendant vertex in $LP_{2, m}$ and $d(v_1, u_2) = 2$, $\partial N_2(v_1)$ = $\{u_2\}$ $\cup$ $\{v_i: d(v_i, v_1) = 2, ~2 \leq i \leq m\}$ = $\{u_2\}$ $\cup$ $\{v_1,v_2,...,v_m\} \setminus \{v_1\}$ and $\partial N_2(u_2)$ = $\{v_j: j$ = 1 to $m\}$ = $\{v_1,v_2,...,v_m\}$. If possible, let $LP_{2,m}$ be 2-DM and $f$ be a 2-DML of $LP_{2,m}$. Then, using the definition of $k$-DML, we get,
		$$ \sum_{u\in N_2(v_1)}{f(u)} =  \sum_{u\in N_2(u_2)}{f(u)}.$$
		$$\Rightarrow  f(u_2) + f(v_2) + f(v_3) + . . . + f(v_m) =  f(v_1) + f(v_2) + . . . + f(v_m).$$	
		$\Rightarrow  f(u_2) =  f(v_1)$  which is a contradiction to the definition of $f$ being a 2-DML of $LP_{2, m}$. Hence $LP_{2,m}$ is not a 2-DM graph in this case.
\item [\rm (ii)]\quad Here, $v_1$ and $v_2$ are pendant vertices in $LP_{n, m}$ and so $d(v_1, v_2) = 2$, $\partial N_2(v_1)$ = $\{u_2\}$ $\cup$ $\{v_i: d(v_i, v_1) = 2, ~2 \leq i \leq m\}$ = $\{u_2\}$ $\cup$ $\{v_1,v_2,...,v_m\} \setminus$ $\{v_1\}$ and $\partial N_2(v_2)$ = $\{u_2\}$ $\cup$ $\{v_1, v_j: d(v_j, v_2) = 2, ~3 \leq j \leq m\}$ = $\{u_2\}$ $\cup$ $\{v_1,v_2,...,v_m\} \setminus \{v_2\}$. If possible, let $LP_{n,m}$ be 2-DM and $f$ be a 2-DML of $LP_{n,m}$. Then, using property (ii) of Theorem \ref{b2}, we get,
	$$ \sum_{v_x\in N[v_1]}{f(v_x)} =  \sum_{v_x\in N[v_2]}{f(v_x)}.$$
	$\Rightarrow  f(v_1) + f(u_1) =  f(v_2) + f(u_1).$	$\Rightarrow  f(v_1) =  f(v_2)$ which is a contradiction to the definition of $f$ being a 2-DML of $LP_{n, m}$. Thus in this case also $LP_{n,m}$ is not a 2-DM graph.
\end{enumerate}
Hence we get the result.
\end{proof}

\begin{cor} \quad \label{b3} {\rm For a given $m \geq 2$ and for every $k\in\mathbb{N}_0$, $LP_{\frac{m(m-1)}{2}+k, m}$ is a $(\frac{m(m-1)}{2}+k)$-DM graph. }
\end{cor}
\begin{proof} \quad This follows from a slight modification of the $n$-DM labeling $f$ of $V(LP_{n, m})$ used in Theorem \ref{b1}. In the modified labeling $f$, vertex $u_{\frac{m(m-1)}{2}+k}$ of $LP_{\frac{m(m-1)}{2}+k, m}$ takes the label $\frac{m(m-1)}{2}$, vertices $v_1,v_2,...,v_m$ take 1,2,...,$m$ and vertices $u_1,u_2,$ ...,$u_{\frac{m(m-1)}{2}+k-1}$ take the remaining elements of $J_{m+\frac{m(m-1)}{2}+k}$. Clearly, it is a $(\frac{m(m-1)}{2}+k)$-DML of $LP_{\frac{m(m-1)}{2}+k, m}$, $k\in\mathbb{N}_0$. Hence we get the result.
\end{proof}

\begin{cor} \quad \label{b4} {\rm For given $m \geq 2$ and $n \geq 3$, the $n$-DM graph $LP_{n, m}$ of least order is $LP_{\frac{m(m-1)}{2}, m}$ where $n$ = $\frac{m(m-1)}{2}$. }
\end{cor}
\begin{proof} \quad Using Theorem \ref{b1}, for $m \geq 2$, and $n \geq 3$, $LP_{n, m}$ is an $n$-DM graph when $m(m-1) \leq 2n.$ The bijective mapping $f$ defined on $V(LP_{n, m})$ in the proof of Theorem \ref{b1} is an $n$-DML of $LP_{n, m}$ and thus the $n$-DML exists on $LP_{n, m}$. Here $u_n$ can not take value less than $\frac{m(m-1)}{2}$ and for the least order $n$-DM graph of $LP_{n,m}$, $u_n$ has to takes the label $\frac{m(m-1)}{2}$, $v_1,v_2,...,v_m$ take labels $1,2,...,m$ and all other labels should be less than $\frac{m(m-1)}{2}$ and greater than $m$. Hence we get the result.
\end{proof}

\begin{rem} \label{b6} One can use Corollary \ref{b3} to obtain different $n$-DMLs of $LP_{n,m}$. Starting with $\frac{m(m-1)}{2}$-DM graph $LP_{\frac{m(m-1)}{2}, m}$, one can produce consecutive super graphs $LP_{\frac{m(m-1)}{2}+h, m}$ which are $(\frac{m(m-1)}{2}+h)$-DM  for $h$ = $1,2,....$ 
\end{rem}

\begin{theorem} \quad \label{b8} {\rm Let $G$ be a graph of order $p \geq 2$ and $u$ be a new vertex. Then $u + G$ is a block if and only if $G$ is connected.}
\end{theorem}
\begin{proof}\quad Let $G$ be connected. Then $G$ has at least one spanning tree, say $T_p$. Clearly, $u+T_p$ is a block. This implies, $u + G$ is a block.

On the otherhand, if $G$ is disconnected, then let $G$ = $G_1$ $\cup$ $G_2$ $\cup$ . . . $\cup$ $G_k$ where $G_1$, $G_2$, . . . , $G_k$ be connected components of $G$ and $k \geq 2$. In this case, $u + G_1$, $u + G_2$, . . . , $u + G_k$ are blocks in $u+G$, $k \geq 2$. And hence $u+G$ is not a block in this case.

Hence we get the result. 
\end{proof}

\begin{cor} \quad \label{b9} {\rm For $m \geq 2$, $LP_{1, m}$ is a block if and only if $<v_1, v_2, . . . , v_m>$ is connected. \hfill $\Box$}
\end{cor}

\section{On $2$-DML of $LP_{1,m}$ = $u_1 + (K_{m-i} \cup K_i)$, $1 \leq i \leq m-i$ and $m \geq 3$}

 In the previous section, we could observe that, in Corollary \ref{b3}, for $k \geq 3$ and $m \geq 2$, once 3-DML of $LP_{3,m}$ is known (its 3-DM labeling mainly depends on the labelings of $v_1, v_2, . . . , v_m$ and of $u_3$), then a $(3+k)$-DML of $LP_{3+k,m}$ can be found easily by considering the same labeling of $v_1, v_2, . . . , v_m$ in $LP_{3+k,m}$ and the labeling of $u_{3+k}$ is same as of $u_3$ in 3-DM of $LP_{3,m}$ and the remaining labels are assigned to the remaining vertices, $u_1, u_2, . . . , u_{2+k}$ for any $k\in\mathbb{N}$. Also, the structure of graph $LP_{n,m}$ depends on the structure of $LP_{1,m}$ whose structure depends on its induced subgraph $<v_1, v_2, . . . , v_m>$ since $LP_{1,m}$ = $u_1$ + $<v_1, v_2, . . . , v_m>$, $m,n\in\mathbb{N}$.  In this section, we concentrate our study on $LP_{1,m}$, $m \geq 3$.

Let $G_i$ be a connected component of $<v_1, v_2, . . . , v_m>$ and $<v_1, v_2, . . . , v_m>$ = $G_1 \cup G_2 \cup . . . \cup G_x$, $1 \leq i \leq x$. Clearly, $u_1 + G_i$ is a block in $LP_{1,m}$, $1 \leq i \leq x$, $m,x\in\mathbb{N}$. 

For $m \geq 2$, $|V(LP_{1,m})|$ = $m+1$ and $dia(LP_{1,m})$ = 1 or 2, $m\in\mathbb{N}$. $dia(LP_{1,m})$ = 1 if and only if $LP_{1,m}$ = $K_{m+1}$. And $LP_{1,m}$ is not DM when $dia(LP_{1,m})$ = 1. 

Let $dia(LP_{1,m})$ = 2, $m \geq 2$ and $m\in\mathbb{N}$. If possible, let $f$ be a DML or 2-DML of $LP_{1,m}$. 

When $LP_{1,m}$ = $P_3$ = $K_1(u_1) + (K_1(v_1) \cup K_1(v_2)) = v_1 u_1 v_2$, $dia(P_3)$ = 2 and $P_3$ is DM but not 2-DM. A path of length 3 and a DML of it are given in Figure 1. 
%%%%%%%Fig.1,2

When $m \geq 2$ and $LP_{1,m}$ contains at least two pendant vertices, $v_i$ and $v_j$, $1 \leq i,j \leq m$, then $dia(LP_{1,m})$  = 2, $LP_{1,m}$ is DM only when $m$ = 2 and is not 2-DM, follows from Theorem \ref{b23}. 

Thus, in the rest of this section, while discussing $2$-DML of $LP_{1,m}$, we consider $dia(LP_{1,m})$  = 2 and $m \geq 3$. Next theorem deals with 2-DML of $LP_{1,m}$ when $x \geq 2$, $<v_1, v_2, . . . , v_m>$ = $K_{m_1} \cup $ $K_{m_2} \cup $ . . . $\cup$ $K_{m_x}$, $1 \leq m_1 \leq m_2 \leq ... \leq m_x$, $m_1+m_2+...+m_x$ = $m$, $m_1+m_2 \geq 3$ (with at the most one pendant vertex in $LP_{1,m}$) and $m_1,m_2,...,m_x,x\in\mathbb{N}$. In this case, each subgraph $K_1(u_1)$ + $K_{m_i}$ = $K_{m_i+1}$ is a block in $LP_{1,m}$, $1 \leq i \leq x$. 

\begin{theorem}\quad \label{d2} {\rm Let $m \geq 3$, $dia(LP_{1,m})$ = 2, $<v_1, v_2, . . . , v_m>$ = $K_{m_1} \cup $ $K_{m_2} \cup $ ... $\cup$ $K_{m_x}$, $x \geq 2$, $1 \leq m_1 \leq m_2 \leq ... \leq m_x$, $m_1+m_2+...+m_x$ = $m$, $m_1+m_2 \geq 3$  and $m_1,m_2,...,m_x,x\in\mathbb{N}$. Then $LP_{1,m}$ is 2-DM if and only if $u_1$ is assigned with a suitable $j$ and $J_{m+1}\setminus \{j\}$ is partitioned into $x$ constant sum partites of orders $m_1,m_2,...,m_x$, $1 \leq j \leq m+1$. }
\end{theorem}
\begin{proof}\quad Let $V(K_{m_i}) = \{ v_{i_j}: j$ = 1 to $m_i\}$, $1 \leq i \leq x$. Then $LP_{1,m}$ is 2-DM if and only if $\exists$ 2-DML $f$ on $LP_{1,m}$ such that 
	$$\sum_{u\in N_2(v_{1_j})}{f(u)} = \sum_{u\in N_2(v_{2_j})}{f(u)} = . . . = \sum_{u\in N_2(v_{x_j})}{f(u)}.$$
	
	$$\Leftrightarrow \sum^{m+1}_{i=1}{i} - (f(u_1) + \sum_{u\in V(K_{m_1})}{f(u)}) = \sum^{m+1}_{i=1}{i} - (f(u_1) + \sum_{u\in V(K_{m_2})}{f(u)})$$
	$$= . . . = \sum^{m+1}_{i=1}{i} - (f(u_1) + \sum_{u\in V(K_{m_x})}{f(u)}).$$ 
	$$\Leftrightarrow \sum_{u\in V(K_{m_1})}{f(u)} =  \sum_{u\in V(K_{m_2})}{f(u)} = . . . = \sum_{u\in V(K_{m_x})}{f(u)}.$$ 
	Thus, under the given conditions, $LP_{1,m}$ is 2-DM if and only if $u_1$ is assigned with a suitable $j$ and $J_{m+1}\setminus \{j\}$ is partitioned into $x$ constant sum partites of orders $m_1,m_2,...,m_x$, $1 \leq j \leq m+1$.
\end{proof}

The necessary condition for bipartition of $J_{m+1} \setminus \{j\}$ with constant sum partites $S_1$ and $S_2$ is that the number of odd numbers in $J_{m+1} \setminus \{j\}$ must be even, $m$ = $m_1+m_2$, $m_1$ = $|S_1|$, $m_2$ = $|S_2|$, $1 \leq j \leq m+1$ and $m \geq 4$. Accordingly, we select $j$ and assign it to $u_1$ so that $LP_{1,m}$ is 2-DM with $m$ = $m_1+m_2$, $m_1$ = $|S_1|$, $m_2$ = $|S_2|$, if it exists, with constant sum $M$ = $\sum_{i\in S_1} i$ = $\sum_{i\in S_2} i$.

Here, by choosing a suitable $j$, we try to obtain possible $x$ constant sum partition of $J_{m+1}\setminus \{j\}$, staring with $x$ = 2,  $1 \leq j \leq m+1$. That is we start with the case of $x$ = 2. In this case, we consider $LP_{1,m}$ such that $dia(LP_{1,m})$ = 2, $<v_1, v_2, . . . , v_m>$ = $K_{m_1}$ $\cup $ $K_{m_2}$, $m_1 \leq m_2$, $m_1+m_2$ = $m \geq 3$ and $m_1,m_2\in\mathbb{N}$. 

In particular, consider $LP_{1,m}$ $\ni$ $m \geq 3$, $dia(LP_{1,m})$ = 2, $<v_1, v_2, . . . , v_m>$ = $K_{1}$ $\cup $ $K_{m-1}$. Figure 3 represents $LP_{1,3}$ with $<v_1, v_2, v_3>$ = $K_{1}$ $\cup $ $K_{2}$ which is 2-DM. The graph $LP_{1,3}$ with its 2-DMLs are given in Figures 4 and 5. 

\vspace{.3cm}	
%Fig.3,4,5
%%%%%%%%%%%%%%%%%%%%%%%%%%%
\begin{center}
	\begin{tikzpicture}
	\node (v1) at (9,2)  [circle,draw,scale=.7]{$v_1$};
	\node (v2) at (10.5,2)  [circle,draw,scale=.7] {$v_2$};
	\node (v3) at (12,2)  [circle,draw,scale=.7] {$v_3$};
	\node (u1) at (10,0)  [circle,draw,scale=.7]{$u_1$};
	
	\draw (u1) -- (v1);
	\draw (u1) -- (v2);
	\draw (u1) -- (v3);
	
	\draw (v2) -- (v3);
	
	%Fig 4	
	\node (y1) at (13,2)  [circle,draw,scale=.7]{3};
	\node (y2) at (14.5,2)  [circle,draw,scale=.7] {1};
	\node (y3) at (16,2)  [circle,draw,scale=.7] {2};
	\node (x1) at (14,0)  [circle,draw,scale=.7]{4};
	
	\draw (x1) -- (y1);
	\draw (x1) -- (y2);
	\draw (x1) -- (y3);
	
	\draw (y2) -- (y3);
	
	%Fig 5	
	\node (z1) at (17,2)  [circle,draw,scale=.7]{4};
	\node (z2) at (18.5,2)  [circle,draw,scale=.7] {1};
	\node (z3) at (20,2)  [circle,draw,scale=.7] {3};
	\node (w1) at (18,0)  [circle,draw,scale=.7]{2};
	
	\draw (w1) -- (z1);
	\draw (w1) -- (z2);
	\draw (w1) -- (z3);
	
	\draw (z2) -- (z3);
	
	\end{tikzpicture}
	
	\vspace{.2cm}		
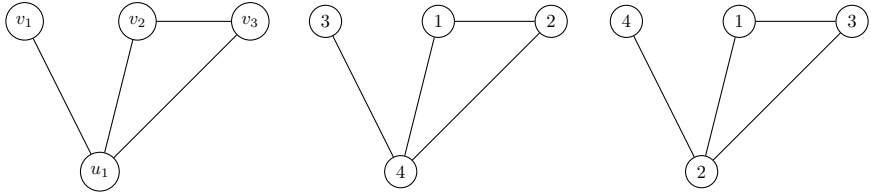
	\footnotesize{ Figure 3. {\small $LP_{1, 3} = K_1 +(K_1 \cup K_2)$} \hspace{.3cm} Figure 4. 2-DML of $LP_{1, 3}$ \hspace{.3cm} Figure 5. 2-DML of $LP_{1, 3}$} 
\end{center}
%%%%%%%%%%%%%%%%%%%%%%%%%%%%%%

\begin{theorem}\quad \label{d3} {\rm Let $m \geq 3$ and $<v_1, v_2, . . . , v_m>$ = $K_{1} \cup $ $K_{m-1}$. Then the graph $LP_{1,m}$ is 2-DM only when $m$ = 3. }
\end{theorem}
\begin{proof} If possible, let $LP_{1,m}$ be 2-DM and $f$ be a 2-DML of it when $<v_1,v_2,...,v_m>$ = $K_{m-1}(v_1, v_2, . . . , v_{m-1})$ $\cup$ $K_{1}(v_m)$ where $K_{n}(v_1,v_2,...,v_{n})$ represents complete graph of order $n$ with the vertex set $\{v_1, v_2, . . . , v_{n}\}$, $m \geq 3$.
	
	Then, for $1 \leq i \leq m-1$,
	$$\sum_{u\in N_2(v_{i})}{f(u)} = \sum_{u\in N_2(v_{m})}{f(u)}.$$
	$$\Rightarrow f(v_m) = \sum^{m-1}_{i=1}{f(v_i)}, ~m \geq 3.$$
	This is possible only when either $\sum^{m-1}_{i=1}{i}$ = $m$ or  $\sum^{m-1}_{i=1}{i}$ = $m+1$ since $v_m$ has to take the biggest number after assigning a suitable value to $u_1$ and the remaining  $m-1$ smaller numbers of $J_{m+1}$ are assigned to vertices  $v_i$, $1 \leq i \leq m-1$. 
	
	This is possible only when $(m-1)m$ = $2m$ or $(m-1)m$ = $2(m+1)$. That is when $m-1 = 2$ or $m^2-3m-2$ = 0. That is when $m = 3$ or $m$ = $\frac{3\pm \sqrt{17}}{2}$. That is when $m = 3$ since $\frac{3\pm \sqrt{17}}{2}\notin\mathbb{N}$.
	
	Hence the result.  
\end{proof}

\begin{theorem}\quad \label{d4} {\rm Let $m \geq 4$ and $<v_1, v_2, . . . , v_m>$ = $K_{2} \cup $ $K_{m-2}$. Then the graph $LP_{1,m}$ is 2-DM only when $m \leq 7$ and it is not 2-DM for $m \geq 8$, $m\in\mathbb{N}$. }
\end{theorem}
\begin{proof}\quad Let us see the result when $m$ = 4, 5, 6, 7.
	
	The necessary condition for bipartition of $J_{m+1} \setminus \{j\}$ with constant sum partites is that the number of odd numbers in $J_{m+1} \setminus \{j\}$ must be even, $1 \leq j \leq m+1$ and $m \geq 4$. Accordingly, we select $j$ and assign it to $u_1$ so that $LP_{1,m}$ is 2-DM, if it exists, with constant sum $M$.
	
	Now, consider constant sum bipartition of $J_{m+1} \setminus \{j\}$ for $m$ = 4,5,6,7,.... In each case of $m$, the number $j$ is identified (starting with biggest possible value of $j$ and then with smaller values) and we represent it by bold font and the remaining numbers are checked for bipartition into constant sum partites $S_1$ and $S_2$ of $J_{m+1} \setminus \{j\}$, $|S_1| \leq |S_2|$, $1 \leq j \leq m+1$.  
	\begin{enumerate}
		\item [\rm 1.] For $m$ = 4 and $<v_1, v_2, . . . , v_m>$ = $K_{2} \cup K_{2}$. 
		\begin{enumerate}
			\item [\rm (i)] $J_{5}$ = $\{1, 2, 3, 4, {\bf 5}\}$. 
			
			Here, $j$ = 5, $S_1$ = $\{1, 4\}$, $S_2$ = $\{2, 3\}$ and $M$ = 5. Corresponding 2-DML of $LP_{1, 4}$ is given in Figure 6. See Figure 6.
			
			\item [\rm (ii)] $J_{5}$ = $\{1, 2, {\bf 3}, 4, 5\}$. 
			
			Here, $j$ = 3, $S_1$ = $\{1, 5\}$, $S_2$ = $\{2, 4\}$ and $M$ = 6. Corresponding 2-DML of $LP_{1, 4}$ is given in Figure 7. See Figure 7.
			
			\item [\rm (iii)] $J_{5}$ = $\{{\bf 1}, 2, 3, 4, 5\}$. 
			
			Here, $j$ = 1, $S_1$ = $\{2, 5\}$, $S_2$ = $\{3, 4\}$ and $M$ = 7. Corresponding 2-DML of $LP_{1, 4}$ is given in Figure 8. See Figure 8.
		\end{enumerate}
		%%%%%%%%%%Fig.6,7,8
		\begin{center}
			\underline{\small{$2$-DMLs of $LP_{1, 4}$ = $u_1 + (K_2 \cup K_2)$ with constant sum $M$}}
			
			\vspace{.2cm}	
			\begin{tikzpicture}
			\node (v1) at (9,2)  [circle,draw,scale=.7]{1};
			\node (v2) at (10,2)  [circle,draw,scale=.7] {4};
			\node (v3) at (11,2)  [circle,draw,scale=.7] {2};
			\node (v4) at (12,2)  [circle,draw,scale=.7] {3};
			\node (u1) at (10.5,0)  [circle,draw,scale=.7]{5};
			
			\draw (u1) -- (v1);
			\draw (u1) -- (v2);
			\draw (u1) -- (v3);
			\draw (u1) -- (v4);
			
			\draw (v1) -- (v2);
			\draw (v3) -- (v4);
			
			%Fig 7	
			\node (y1) at (13,2)  [circle,draw,scale=.7]{1};
			\node (y2) at (14,2)  [circle,draw,scale=.7] {5};
			\node (y3) at (15,2)  [circle,draw,scale=.7] {2};
			\node (y4) at (16,2)  [circle,draw,scale=.7] {4};
			\node (x1) at (14.5,0)  [circle,draw,scale=.7]{3};
			
			\draw (x1) -- (y1);
			\draw (x1) -- (y2);
			\draw (x1) -- (y3);
			\draw (x1) -- (y4);
			
			\draw (y1) -- (y2);
			\draw (y3) -- (y4);
			
			%Fig 8	
			\node (z1) at (17,2)  [circle,draw,scale=.7]{2};
			\node (z2) at (18,2)  [circle,draw,scale=.7] {5};
			\node (z3) at (19,2)  [circle,draw,scale=.7] {3};
			\node (z4) at (20,2)  [circle,draw,scale=.7] {4};
			\node (w1) at (18.5,0)  [circle,draw,scale=.7]{1};
			
			\draw (w1) -- (z1);
			\draw (w1) -- (z2);
			\draw (w1) -- (z3);
			\draw (w1) -- (z4);
			
			\draw (z1) -- (z2);
			\draw (z3) -- (z4);
			
			\end{tikzpicture}
			
			\vspace{.2cm}		
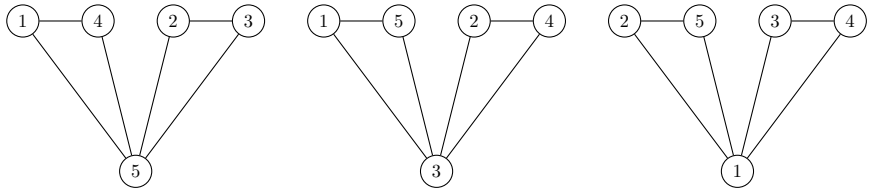
			\footnotesize{\hspace{-.4cm} Fig. 6. 2-DML with $M$ = 5. \hspace{.1cm} Fig. 7. 2-DML with $M$ = 6. \hspace{.1cm} Fig. 8. 2-DML with $M$ = 7. }
		\end{center}
		%%%%%%%%%%%%%%%%%%%%%%%%%%%%%%
		\item [\rm 2.] For $m$ = 5 and $<v_1, v_2, . . . , v_m>$ = $K_{2} \cup K_{3}$. 
		\begin{enumerate}
			\item [\rm (i)] $J_{6}$ = $\{1, 2, 3, 4, {\bf 5}, 6\}$. 
			
			Here, $j$ = 5, $S_1$ = $\{2, 6\}$, $S_2$ = $\{1, 3, 4\}$ and $M$ = 8. Corresponding 2-DML of $LP_{1, 5}$ is given in Figure 9. See Figure 9.
			
			\item [\rm (ii)] $J_{6}$ = $\{1, 2, {\bf 3}, 4, 5, 6\}$. 
			
			Here, $j$ = 3, $S_1$ = $\{4, 5\}$, $S_2$ = $\{1, 2, 6\}$ and $M$ = 9. Corresponding 2-DML of $LP_{1, 5}$ is given in Figure 10. See Figure 10.
			
			\item [\rm (iii)] $J_{6}$ = $\{{\bf 1}, 2, 3, 4, 5, 6\}$. 
			
			Here, $j$ = 1, $S_1$ = $\{4, 6\}$, $S_2$ = $\{2, 3, 5\}$ and $M$ = 10. Corresponding 2-DML of $LP_{1, 5}$ is given in Figure 11. See Figure 11.
		\end{enumerate}
		%%%%%%%%%%%Fig.9,10,11
		\begin{center}
			\underline{\small{$2$-DMLs of $LP_{1, 5}$ = $u_1 + (K_3 \cup K_2)$ with constant sum $M$}}
			
			\vspace{.2cm}	
			\begin{tikzpicture}
			\node (v1) at (10,2)  [circle,draw,scale=.7]{2};
			\node (v2) at (11,2)  [circle,draw,scale=.7] {6};
			\node (v3) at (12,2)  [circle,draw,scale=.7] {1};
			\node (v4) at (12,1)  [circle,draw,scale=.7] {3};
			\node (v5) at (13,1)  [circle,draw,scale=.7] {4};
			\node (u1) at (11,0)  [circle,draw,scale=.7]{5};
			
			\draw (u1) -- (v1);
			\draw (u1) -- (v2);
			\draw (u1) -- (v3);
			\draw (u1) -- (v4);
			\draw (u1) -- (v5);
			
			\draw (v1) -- (v2);
			\draw (v3) -- (v4);
			\draw (v3) -- (v5);
			\draw (v4) -- (v5);
			
			%Fig 10	
			\node (y1) at (14,2)  [circle,draw,scale=.7]{4};
			\node (y2) at (15,2)  [circle,draw,scale=.7] {5};
			\node (y3) at (16,2)  [circle,draw,scale=.7] {1};
			\node (y4) at (16,1)  [circle,draw,scale=.7] {2};
			\node (y5) at (17,1)  [circle,draw,scale=.7] {6};
			\node (x1) at (15,0)  [circle,draw,scale=.7]{3};
			
			\draw (x1) -- (y1);
			\draw (x1) -- (y2);
			\draw (x1) -- (y3);
			\draw (x1) -- (y4);
			\draw (x1) -- (y5);
			
			\draw (y1) -- (y2);
			\draw (y3) -- (y4);
			\draw (y3) -- (y5);
			\draw (y4) -- (y5);
			
			%Fig 11	
			\node (z1) at (18,2)  [circle,draw,scale=.7]{4};
			\node (z2) at (19,2)  [circle,draw,scale=.7] {6};
			\node (z3) at (20,2)  [circle,draw,scale=.7] {2};
			\node (z4) at (20,1)  [circle,draw,scale=.7] {3};
			\node (z5) at (20.8,1)  [circle,draw,scale=.7] {5};
			\node (w1) at (19,0)  [circle,draw,scale=.7]{1};
			
			\draw (w1) -- (z1);
			\draw (w1) -- (z2);
			\draw (w1) -- (z3);
			\draw (w1) -- (z4);
			\draw (w1) -- (z5);

			\draw (z1) -- (z2);
			\draw (z3) -- (z4);
			\draw (z3) -- (z5);
			\draw (z4) -- (z5);
			
			\end{tikzpicture}
			
			\vspace{.2cm}		
			\footnotesize{ Figure 9. $M$ = 8. \hspace{1cm} Figure 10.  $M$ = 9. \hspace{1cm} Figure 11. $M$ = 10. }
		\end{center}
		%%%%%%%%%%%%%%%%%%%%%%%%%%%%%%
		\item [\rm 3.] For $m$ = 6 and $<v_1, v_2, . . . , v_m>$ = $K_{2} \cup K_{4}$. 
		\begin{enumerate}
			\item [\rm (i)] $J_{7}$ = $\{1, 2, 3, 4, 5, {\bf 6}, 7\}$. 
			
			Here, $j$ = 6, $S_1$ = $\{4, 7\}$, $S_2$ = $\{1, 2, 3, 5\}$ and $M$ = 11. Corresponding 2-DML of $LP_{1, 6}$ is given in Figure 12. See Figure 12.
			
			\item [\rm (ii)] $J_{7}$ = $\{1, 2, 3, {\bf 4}, 5, 6, 7\}$. 
			
			Here, $j$ = 4, $S_1$ = $\{5, 7\}$, $S_2$ = $\{1, 2, 3, 6\}$ and $M$ = 12. Corresponding 2-DML of $LP_{1, 6}$ is given in Figure 13. See Figure 13.
			
			\item [\rm (iii)] $J_{7}$ = $\{1, {\bf 2}, 3, 4, 5, 6, 7\}$. 
			
			Here, $j$ = 2, $S_1$ = $\{6, 7\}$, $S_2$ = $\{1, 3, 4, 5\}$ and $M$ = 13. Corresponding 2-DML of $LP_{1, 6}$ is given in Figure 14. See Figure 14.
		\end{enumerate}
		%%%%%%%%%Fig.12,13,14
		\begin{center}
			\underline{\small{$2$-DMLs of $LP_{1, 6}$ = $u_1 + (K_4 \cup K_2)$ with constant sum $M$}}
			
			\vspace{.2cm}	
			\begin{tikzpicture}
			\node (v1) at (10,2)  [circle,draw,scale=.7]{4};
			\node (v2) at (11,2)  [circle,draw,scale=.7] {7};
			\node (v3) at (11.75,2)  [circle,draw,scale=.7] {1};
			\node (v4) at (12,1)  [circle,draw,scale=.7] {2};
			\node (v5) at (13,1.5)  [circle,draw,scale=.7] {3};
			\node (v6) at (13,.5)  [circle,draw,scale=.7] {5};
			\node (u1) at (11.5,0)  [circle,draw,scale=.7]{6};
			
			\draw (u1) -- (v1);	\draw (u1) -- (v2);
			\draw (u1) -- (v3);	\draw (u1) -- (v4);
			\draw (u1) -- (v5);	\draw (u1) -- (v6);
			
			\draw (v1) -- (v2);	\draw (v3) -- (v4);
			\draw (v3) -- (v5);	\draw (v3) -- (v6);
			\draw (v4) -- (v5);	\draw (v4) -- (v6);
			\draw (v5) -- (v6);
			
			%Fig 10	
			\node (y1) at (14,2)  [circle,draw,scale=.7]{5};
			\node (y2) at (15,2)  [circle,draw,scale=.7] {7};
			\node (y3) at (15.75,2)  [circle,draw,scale=.7] {1};
			\node (y4) at (15.75,1)  [circle,draw,scale=.7] {2};
			\node (y5) at (17,1.5)  [circle,draw,scale=.7] {3};
			\node (y6) at (17,.5)  [circle,draw,scale=.7] {6};
			\node (x1) at (15,0)  [circle,draw,scale=.7]{4};
			
			\draw (x1) -- (y1);	\draw (x1) -- (y2);
			\draw (x1) -- (y3);	\draw (x1) -- (y4);
			\draw (x1) -- (y5);	\draw (x1) -- (y6);
			
			\draw (y1) -- (y2);
			\draw (y3) -- (y4);	\draw (y3) -- (y5);
			\draw (y3) -- (y6);	\draw (y4) -- (y5);
			\draw (y4) -- (y6);  \draw (y5) -- (y6);
			
			%Fig 11	
			\node (z1) at (18,2)  [circle,draw,scale=.7]{6};
			\node (z2) at (19,2)  [circle,draw,scale=.7] {7};
			\node (z3) at (19.75,2)  [circle,draw,scale=.7] {1};
			\node (z4) at (19.75,1)  [circle,draw,scale=.7] {3};
			\node (z5) at (21,1.5)  [circle,draw,scale=.7] {4};
			\node (z6) at (21,.5)  [circle,draw,scale=.7] {5};
			\node (w1) at (19,0)  [circle,draw,scale=.7]{2};
			
			\draw (w1) -- (z1);
			\draw (w1) -- (z2);
			\draw (w1) -- (z3);
			\draw (w1) -- (z4);
			\draw (w1) -- (z5);
			\draw (w1) -- (z6);

			\draw (z1) -- (z2);
			\draw (z3) -- (z4);
			\draw (z3) -- (z5);
			\draw (z3) -- (z6);
			\draw (z4) -- (z5);
			\draw (z4) -- (z6);
			\draw (z5) -- (z6);
			
			\end{tikzpicture}
			
			\vspace{.2cm}		
			\footnotesize{ Figure 12. $M$ = 11. \hspace{1cm} Figure 13.  $M$ = 12. \hspace{1cm} Figure 14. $M$ = 13. }
			
		\end{center}
	%%%%%%%%%%%%%%%%%%%
		\item [\rm 4.] For $m$ = 7 and $<v_1, v_2, . . . , v_m>$ = $K_{2} \cup K_{5}$. 
		
		$J_{8}$ = $\{1, 2, 3, 4, 5, {\bf 6}, 7, 8\}$ and $S_1$ = $\{7, 8\}$ and $S_2$ = $\{1, 2, 3, 4, 5\}$ are its constant sum bipartite subsets with $M$ = 15. The corresponding 2-DML of $LP_{1, 7}$ is given in Figure 15.  See Figure 15. 
		%%%%%%%%%%%%Fig.15
		\vspace{.2cm}		
		\begin{center}	
			\begin{tikzpicture}
			\node (v1) at (10.5,1.5)  [circle,draw,scale=.7]{8};
			\node (v2) at (11.5,1.5)  [circle,draw,scale=.7] {7};
			\node (v3) at (12.75,1.5)  [circle,draw,scale=.7] {1};
			\node (v4) at (13.75,.75)  [circle,draw,scale=.7] {2};
			\node (v5) at (14.25,1.5)  [circle,draw,scale=.7] {3};
			\node (v6) at (15,.75)  [circle,draw,scale=.7] {4};
			\node (v7) at (13.85,0)  [circle,draw,scale=.7] {5};
			\node (u1) at (11.5,-.5)  [circle,draw,scale=.7]{6};
			
			\draw (u1) -- (v1);	\draw (u1) -- (v2);
			\draw (u1) -- (v3);	\draw (u1) -- (v4);
			\draw (u1) -- (v5);	\draw (u1) -- (v6);
			\draw (u1) -- (v7);
			
			\draw (v1) -- (v2);	\draw (v3) -- (v4);
			\draw (v3) -- (v5);	\draw (v3) -- (v6);
			\draw (v3) -- (v7); \draw (v4) -- (v5);
			\draw (v4) -- (v6); \draw (v4) -- (v7);
			\draw (v5) -- (v6); \draw (v5) -- (v7);
			\draw (v6) -- (v7);
			
			\end{tikzpicture}
			
			\vspace{.2cm}		
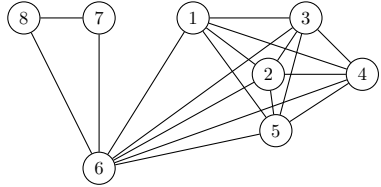
			\footnotesize{Figure 15. $2$-DMLs of $LP_{1, 7}$ = $u_1 + (K_5 \cup K_2)$ with $M$ = 15. }
		\end{center}
		%%%%%%%%%%%%%%%%%%%%%%
		\item [\rm 5.] For $m$ = 8 and $<v_1, v_2, . . . , v_m>$ = $K_{2} \cup K_{6}$. 
		
		$J_{9}$ = $\{1, 2, 3, 4, 5, 6, 7, 8, 9\}$ and the number of odd numbers in $J_9$ is 5, an odd number and so $u_1$ takes an odd number, say, $j$ for possible constant sum bipartition of $J_9 \setminus \{j\}$. Under these conditions, sum of the two biggest numbers 9+8 is less than sum of the six smallest numbers, 1+2+...+6 = 21 and hence 2-DML is not possible in this case. This is also true for higher values of $m$.
	\end{enumerate}
	
	Thus from the above cases, we get, for $m-2 \geq 2$, $LP_{1,m}$ = $u_1$ + $(K_{2} \cup K_{m-2})$ is 2-DM only for $m$ = 4 to 7 and is not 2-DM when $m \geq 8$, $m\in\mathbb{N}$. Hence we get the result.  
\end{proof}

\begin{theorem}\quad \label{d5} {\rm Let $m \geq 6$ and $<v_1, v_2, . . . , v_m>$ = $K_{3} \cup $ $K_{m-3}$. Then the graph $LP_{1,m}$ is 2-DM only when $m \leq 10$ and not 2-DM for $m \geq 11$, $m\in\mathbb{N}$. }
\end{theorem}
\begin{proof}\quad Let us consider $LP_{1,m}$ for $m \geq 6$ with $<v_1, v_2, . . . , v_m>$ = $K_{3} \cup K_{m-3}$. Here, we start with $m-3$ = 3, the smallest possible value of $m-3$ when $m-3 \geq 3$. The bold font number in each $J_{m+1}$ indicates that it is the possible label $j$ for $u_1$ in $LP_{1, m}$ to obtain constant sum partites $S_1$ and $S_2$ of $J_{m+1} \setminus \{j\}$, $|S_1| \leq |S_2|$ and $M$ is the constant sum.  
	\begin{enumerate}
		\item [\rm 1.] For $m$ = 6 and $<v_1, v_2, . . . , v_m>$ = $K_{3} \cup K_{3}$. 
		\begin{enumerate}
			\item [\rm (i)] $J_{7}$ = $\{1, 2, 3, 4, 5, {\bf 6}, 7\}$. 
			
			Here, $j$ = 6, $S_1$ = $\{1, 3, 7\}$,  $S_2$ = $\{2, 4, 5\}$ and $M$ = 11.
			
			\item [\rm (ii)] $J_{7}$ = $\{1, 2, 3, {\bf 4}, 5, 6, 7\}$. 
			
			Here, $j$ = 4, $S_1$ = $\{1, 5, 6\}$, $S_2$ = $\{2, 3, 7\}$ and $M$ = 12.
			
			\item [\rm (iii)] $J_{7}$ = $\{1, {\bf 2}, 3, 4, 5, 6, 7\}$. 
			
			Here, $j$ = 2, $S_1$ = $\{1, 5, 7\}$, $S_2$ = $\{3, 4, 6\}$ and $M$ = 13.
		\end{enumerate}
		
		\item [\rm 2.] For $m$ = 7 and $<v_1, v_2, . . . , v_m>$ = $K_{3} \cup K_{4}$. 
		\begin{enumerate}
			\item [\rm (i)] $J_{8}$ = $\{1, 2, 3, 4, 5, 6, 7, {\bf 8}\}$. 
			
			Here, $j$ = 8, $M$ = 14 and possible values of $S_1$ and $S_2$ are
			
			$S_1$ = $\{1, 6, 7\}$, $S_2$ = $\{2, 3, 4, 5\}$; 
			
			$S_1$ = $\{2, 5, 7\}$, $S_2$ = $\{1, 3, 4, 6\}$;
			
			$S_1$ = $\{3, 5, 6\}$, $S_2$ = $\{1, 2, 4, 7\}$.
			
			\item [\rm (ii)] $J_{8}$ = $\{1, 2, 3, 4, 5, {\bf 6}, 7, 8\}$. 
			
			Here, $j$ = 6, $M$ = 15 and possible values of $S_1$ and $S_2$ are
			
			$S_1$ = $\{2, 5, 8\}$, $S_2$ = $\{1, 3, 4, 7\}$;  
			
			$S_1$ = $\{3, 4, 8\}$, $S_2$ = $\{1, 2, 5, 7\}$;
			
			$S_1$ = $\{3, 5, 7\}$, $S_2$ = $\{1, 2, 4, 8\}$.
			
			\item [\rm (iii)] $J_{8}$ = $\{1, 2, 3, {\bf 4}, 5, 6, 7, 8\}$. 
			
			Here, $j$ = 4, $M$ = 16 and possible values of $S_1$ and $S_2$ are
			
			$S_1$ = $\{1, 7, 8\}$, $S_2$ = $\{2, 3, 5, 6\}$;  
			
			$S_1$ = $\{2, 6, 8\}$, $S_2$ = $\{1, 3, 5, 7\}$;
			
			$S_1$ = $\{3, 6, 7\}$, $S_2$ = $\{1, 2, 5, 8\}$.
			
			\item [\rm (iv)] $J_{8}$ = $\{1, {\bf 2}, 3, 4, 5, 6, 7, 8\}$. 
			
			Here, $j$ = 2, $M$ = 17 and possible values of $S_1$ and $S_2$ are
			
			$S_1$ = $\{4, 6, 7\}$, $S_2$ = $\{1, 3, 5, 8\}$;  
			
			$S_1$ = $\{3, 6, 8\}$, $S_2$ = $\{1, 4, 5, 7\}$;
			
			$S_1$ = $\{4, 5, 8\}$, $S_2$ = $\{1, 3, 6, 7\}$.
		\end{enumerate}
		
		\item [\rm 3.] For $m$ = 8 and $<v_1, v_2, . . . , v_m>$ = $K_{3} \cup K_{5}$. 
		\begin{enumerate}
			\item [\rm (i)] $J_{9}$ = $\{1, 2, 3, 4, 5, 6, 7, 8, {\bf 9}\}$. 
			
			Here, $j$ = 9, $M$ = 18 and possible values of $S_1$ and $S_2$ are
			
			$S_1$ = $\{3, 7, 8\}$, $S_2$ = $\{1, 2, 4, 5, 6\}$;  
			
			$S_1$ = $\{4, 6, 8\}$, $S_2$ = $\{1, 2, 3, 5, 7\}$;
			
			$S_1$ = $\{5, 6, 7\}$, $S_2$ = $\{1, 2, 3, 4, 8\}$.
			
			\item [\rm (ii)] $J_{9}$ = $\{1, 2, 3, 4, 5, 6, {\bf 7}, 8, 9\}$. 
			
			Here, $j$ = 7, $M$ = 19 and possible values of $S_1$ and $S_2$ are
			
			$S_1$ = $\{2, 8, 9\}$, $S_2$ = $\{1, 3, 4, 5, 6\}$;  
			
			$S_1$ = $\{4, 6, 9\}$, $S_2$ = $\{1, 2, 3, 5, 8\}$;
			
			$S_1$ = $\{5, 6, 8\}$, $S_2$ = $\{1, 2, 3, 4, 9\}$.
			
			\item [\rm (iii)] $J_{9}$ = $\{1, 2, 3, 4, {\bf 5}, 6, 7, 8, 9\}$. 
			
			Here, $j$ = 5, $M$ = 20 and possible values of $S_1$ and $S_2$ are 
			
			$S_1$ = $\{3, 8, 9\}$, $S_2$ = $\{1, 2, 4, 6, 7\}$;  
			
			$S_1$ = $\{4, 7, 9\}$, $S_2$ = $\{1, 2, 3, 6, 8\}$.
			
			\item [\rm (iv)] $J_{9}$ = $\{1, 2, {\bf 3}, 4, 5, 6, 7, 8, 9\}$. 
			
			Here, $j$ = 3, $M$ = 21 and possible values of $S_1$ and $S_2$ are 
			
			$S_1$ = $\{4, 8, 9\}$, $S_2$ = $\{1, 2, 5, 6, 7\}$;  
			
			$S_1$ = $\{5, 7, 9\}$, $S_2$ = $\{1, 2, 4, 6, 8\}$;
			
			$S_1$ = $\{6, 7, 8\}$, $S_2$ = $\{1, 2, 4, 5, 9\}$.
			
			\item [\rm (iv)] $J_{9}$ = $\{{\bf 1}, 2, 3, 4, 5, 6, 7, 8, 9\}$. 
			
			Here, $j$ = 1, $M$ = 22 and possible values of $S_1$ and $S_2$ are  
			
			$S_1$ = $\{5, 8, 9\}$, $S_2$ = $\{2, 3, 4, 6, 7\}$; 
			
			$S_1$ = $\{6, 7, 9\}$, $S_2$ = $\{2, 3, 4, 5, 8\}$.
		\end{enumerate}
		
		\item [\rm 4.] For $m$ = 9 and $<v_1, v_2, . . . , v_m>$ = $K_{3} \cup K_{6}$. 
		\begin{enumerate}
			\item [\rm (i)] $J_{10}$ = $\{1, 2, 3, 4, 5, 6, 7, 8, {\bf 9}, 10\}$. 
			
			Here, $j$ = 9, $M$ = 23 and possible values of $S_1$ and $S_2$ are 
			
			$S_1$ = $\{5, 8, 10\}$, $S_2$ = $\{1, 2, 3, 4, 6, 7\}$;  
			
			$S_1$ = $\{6, 7, 10\}$, $S_2$ = $\{1, 2, 3, 4, 5, 8\}$.
			
			\item [\rm (ii)] $J_{10}$ = $\{1, 2, 3, 4, 5, 6, {\bf 7}, 8, 9, 10\}$. 
			
			Here, $j$ = 7, $M$ = 24 and possible values of $S_1$ and $S_2$ are 
			
			$S_1$ = $\{5, 9, 10\}$, $S_2$ = $\{1, 2, 3, 4, 6, 8\}$;  
			
			$S_1$ = $\{6, 8, 10\}$, $S_2$ = $\{1, 2, 3, 4, 5, 9\}$.
			
			\item [\rm (iii)] $J_{10}$ = $\{1, 2, 3, 4, {\bf 5}, 6, 7, 8, 9, 10\}$. 
			
			Here, $j$ = 5, $M$ = 25 and possible values of $S_1$ and $S_2$ are 
			
			$S_1$ = $\{6, 9, 10\}$, $S_2$ = $\{1, 2, 3, 4, 7, 8\}$;  
			
			$S_1$ = $\{7, 8, 10\}$, $S_2$ = $\{1, 2, 3, 4, 6, 9\}$.
			
			\item [\rm (iv)] $J_{10}$ = $\{1, 2, {\bf 3}, 4, 5, 6, 7, 8, 9, 10\}$. 
			
			Here, $j$ = 3, $S_1$ = $\{7, 9, 10\}$, $S_2$ = $\{1, 2, 4, 5, 6, 8\}$ and $M$ = 26.
			
			\item [\rm (iv)] $J_{10}$ = $\{{\bf 1}, 2, 3, 4, 5, 6, 7, 8, 9, 10\}$. 
			
			Here, $j$ = 1, $S_1$ = $\{8, 9, 10\}$, $S_2$ = $\{2, 3, 4, 5, 6, 7\}$ and $M$ = 27.
		\end{enumerate}
		
		\item [\rm 5.] For $m$ = 10 and $<v_1, v_2, . . . , v_m>$ = $K_{3} \cup K_{7}$. 
		\begin{enumerate}
			\item [\rm (i)] $J_{11}$ = $\{1, 2, 3, 4, 5, 6, 7, 8, 9, {\bf 10}, 11\}$. 
			
			Here, $j$ = 10, $S_1$ = $\{8, 9, 11\}$, $S_2$ = $\{1, 2, 3, 4, 5, 6, 7\}$ and $M$ = 28.
			
			\item [\rm (ii)] $J_{11}$ = $\{1, 2, 3, 4, 5, 6, 7, {\bf 8}, 9, 10, 11\}$. 
			
			In this case, $\frac{S}{2} = \frac{1}{2}(\sum^{11}_{i=1}{i}-8) = 29$ and sum of any 3 distinct numbers of $J_{11} \setminus \{8\}$ $\neq$ 29 and so constant sum bipartition doesn't exist and thereby $S_1$ and $S_2$ doesn't exist. 
			
			\item [\rm (iii)] $J_{11}$ = $\{1, 2, 3, 4, 5, {\bf 6}, 7, 8, 9, 10, 11\}$. 
			
			Here, $j$ = 6, $S_1$ = $\{9, 10, 11\}$, $S_2$ = $\{1, 2, 3, 4, 5, 7, 8\}$ and $M$ = 30.
			
			\item [\rm (iv)] $J_{11}$ = $\{1, 2, 3, {\bf 4}, 5, 6, 7, 8, 9, 10, 11\}$. 
			
			In this case, sum of 7 smallest numbers = 32 which is greater than 30 = 11+10+9, sum of 3 biggest numbers in $J_{11}\setminus \{4\}$ and hence $S_1$ and $S_2$ doesn't exist.
		\end{enumerate}
		\item [\rm 6.] For $m$ = 11 and $<v_1, v_2, . . . , v_m>$ = $K_{3} \cup K_{8}$. 
		\begin{enumerate}
			\item [\rm (i)] $J_{12}$ = $\{1, 2, 3, 4, 5, 6, 7, 8, 9, 10, 11, {\bf 12}\}$. 
			
			In this case, $\frac{S}{2} = \frac{1}{2}(\sum^{12}_{i=1}{i}-12) = 33$ and sum of 8 smallest numbers of $J_{11}$ = 36 which is greater than 30 = 11+10+9, sum of 3 biggest numbers in $J_{11}$ and hence $S_1$ and $S_2$ doesn't exist.
		\end{enumerate}	
		The above also implies, when $m$ is greater than 11, sum of $m-3$ smallest numbers of $J_{m+1}$ will be greater than sum of 3 biggest numbers and thereby constant sum bipartition doesn't exist for $J_{m+1}\setminus \{j\}$, $1 \leq j \leq m+1$.
	\end{enumerate}
	Hence we get the result.
\end{proof}

 For $m_1$ = $|S_1|$ = 1, $m$ = 3 is the only possible value of $m$ for which $LP_{1, 3}$ = $u_1 + (K_1(v_1) \cup K_2(v_2v_3))$ is 2-DM. The  graph and its two 2-DMLs are given in Figures 3, 4, 5.

For $m_1$ = 2, $m$ = 4,5,6,7 are the possible values of $m$ for which $LP_{1, m}$ = $u_1 + (K_2 \cup K_{m-2})$ is 2-DM and the corresponding 2-DM labeled graphs are given in Figures 6 to 14.

For $m_1$ = 3, $m$ = 6,7,8,9,10 are the possible values of $m$ for which $LP_{1, m}$ = $u_1 + (K_3$ $\cup$ $K_{m-3})$ is 2-DM. In each case, we provide constant sum bipartition sets.  

For $m_1$ = 4 to 22, we calculate possible values of $m$ for which $LP_{1, m}$ = $u_1 + (K_4 \cup K_{m-4})$ is 2-DM and these values are presented in Table 1. Detailed calculations for obtaining these values of $m$ are given in the Annexure.

\begin{theorem}\quad \label{d6} {\rm For $n \geq 3$, $J_n$ is constant sum bipartite if and only if $n \equiv ~ 0,3 ~(mod ~4)$ if and only if $J_n$ contains even number of odd numbers. }
\end{theorem} 
\begin{proof}\quad Consider the case, $n \equiv ~ 0 ~(mod ~4)$. Let $n$ = $4m$ and $A_i$ = $\{i, 4m+1-i\}$ for $i$ = 1 to $2m$. The two elements of $A_i$ are odd and even, $i$ = 1 to $2m$ and $m\in\mathbb{N}$. Take union of any $m$ number of $A_i$ sets as $S_1$ and union of the remaining $m$ sets as $S_2$. Clearly, $S_1$ and $S_2$ partition the set $J_{4m}$ into constant sum bipartition with constant sum $M$ = $m(4m+1)$. Moreover, the number of odd numbers in $J_{4m}$ is $2m$, an even number. 
	
 Consider the case, $n \equiv ~ 3 ~(mod ~4)$.  Let $n$ = $4m+3$, $A_i$ = $\{i, 4m+3-i\}$ for $i$ = 1 to $2m+1$, and $A_{2m+2}$ = $\{4m+3\}$, $m\in\mathbb{N}$. Here, there are $2m+2$ number of $A_i$s such that sum of the elements in each $A_i$ is a constant equal to $4m+3$, $1 \leq i \leq 2m+2$.  Take $S_1$ as union of any $m+1$ sets of $A_i$s and union of the remaining $m+1$ sets as $S_2$. Then, $S_1$ and $S_2$ are constant sum bipartites of $J_{4m+3}$ with constant sum $M$ = $(m+1)(4m+3)$. Moreover, the number of odd numbers in $J_{4m+3}$ is $2m+2$, an even number. 
	 
	 Thus, in these two cases constant sum bipartition exists on $J_n$ and the number of odd numbers in each case is even. 
	 
	 On the otherhand, when $n \equiv ~ 1,2 ~(mod ~4)$, the number of odd numbers in $J_n$ = $J_{4m+1}$ as well as in $J_{4m+2}$ is $2m+1$ which is an odd number and so the elements of $J_{4m+1}$ as well as of $J_n$ = $J_{4m+2}$ can not be partitioned into $S_1$ and $S_2$ such that sum of the elements in $S_1$ is same as the sum of the elements in $S_2$. Thus, in these two cases constant sum bipartition of $J_{n}$ is not possible.
	 
	 Combining the above cases, we get the result.
\end{proof}
%\begin{note}\quad \label{d7} {\rm For $n \geq 3$ and $1 \leq j \leq n+1$, $J_{n+1}\setminus \{j\}$ is constant sum bipartite for some $j$ if and only if $J_{n+1}\setminus \{j\}$ contains even number of odd numbers. }\end{note} 

\vspace{.1cm}
\noindent
{\bf To find $m$ when $m_1$ is given in 2-DM graph $LP_{1, m} = K_1(u_1)+(K_{m_1}\cup K_{m_2})$, 
	
	\hfill $m$ = $m_1 + m_2 \geq 3$ and $1 \leq m_1 \leq m_2$}.

\vspace{.1cm}
Theorems \ref{d2} and \ref{d6} ensures existence of 2-DM graphs $LP_{1, m}$ when $LP_{1, m}$ = $K_1(u_1)+(K_{m_1}\cup K_{m_2})$ but not for all possible values of $m_1$ and $m_2$, $m$ = $m_1 + m_2 \geq 3$ and $1 \leq m_1 \leq m_2$. Here, we try to find out the possible values of $m_1$ = $|S_1|$ and $m_2$ = $|S_2|$ for which $LP_{1, m}$ = $K_1(u_1)+(K_{m_1}\cup K_{m_2})$ is 2-DM. Let $f$ be a 2-DML of $LP_{1,m} = u_1+(K_{m_1}\cup K_{m_2})$, $m$ = $m_1 + m_2 \geq 3$ and $1 \leq m_1 \leq m_2$. Let $j$ = $f(u_1)\in J_{m+1}$ be a suitable value such that $J_{m+1} \setminus \{j\}$ can be bipartitioned into constant sum partites $S_1$ and $S_2$, $m_1 = |S_1|$ and $m_2 = |S_2|$, $m_1 \leq m_2$ and $1 \leq j \leq m+1$. Such $S_1$ and $S_2$ exist by Theorems \ref{d2} and \ref{d6}. Here, we obtain $S_1$ and $S_2$ of $J_{m+1} \setminus \{j\}$, starting with the smallest possible $m$ and assigning $f(u_1)$ = $m+1$ = $j$ when $m+1 \equiv ~0,1~ (mod~ 4)$ or $f(u_1)$ = $m$ = $j$ when $m+1 \equiv ~2,3~ (mod~ 4)$, $m \geq 3$. For $m_1$ = 1 to 22, Table 1 presents possible such values of $m$ for which $LP_{1,m} = u_1+(K_{m_1}\cup K_{m_2})$ is 2-DM, $m$ = $m_1+m_2$ and $1 \leq m_1 \leq m_2$. Detailed calculations for obtaining these values of $m$ (which are presented in Table 1) are given in the Annexure. From these values of $S_1$ and $S_2$ of $J_{m+1} \setminus \{j\}$, we obtain 2-DM graphs $LP_{1, m} = u_1+(K_{m_1}\cup K_{m_2})$ by assigning the elements of $S_1$ as the vertex labels of $K_{m_1}$, the elements of $S_2$ as the vertex labels of $K_{m_2}$ and $f(u)$ = $j$, $m$ = $m_1 + m_2 \geq 3$ and $1 \leq m_1 \leq m_2$. It is also noted that a general formula to obtain such values of $m$ seems to be difficult. 
%%%%%%Table-1

\begin{oprm}\quad \label{p1} {\rm For a given value of $m_1$, find a general formula for $m$ for which $LP_{1, m} = u_1+(K_{m_1}\cup K_{m_2})$ is 2-DM, $1 \leq m_1 \leq m_2$, $m = m_1 + m_2 \geq 3$ and $m,m_1,m_2\in\mathbb{N}$. }
\end{oprm} 

\begin{table} \label{t1} 
	\caption{\small{Value(s) of $m$ when $LP_{1, m}$ = $u_1+(K_{m_1} \cup K_{m_2})$ and is 2-DM.}}
	\begin{center}
		%\scalebox{0.5}{
		\begin{tabular}{|c|c|c|} \hline \hline
			\hspace{.1cm} $m_1$ \hspace{.1cm}
			& $m_2$ = $m-m_1 \geq m_1$  & \hspace{.3cm} Possible value(s) of $m$ \hspace{.3cm}  \\\hline \hline 
			1 & $m-1$  & $m$ = 3   \\\hline\hline 
			2 & $m-2$  & $m$ = 4 to 7 \\\hline
			3 & $m-3$  & $m$ = 6 to 10 \\\hline
			4 & $m-4$  & $m$ = 8 to 13 \\\hline
			5 & $m-5$  & $m$ = 10 to 16 \\\hline\hline
			6 & $m-6$  & $m$ = 12 to 20 \\\hline
			7 & $m-7$  & $m$ = 14 to 23 \\\hline
			8 & $m-8$  & $m$ = 16 to 26 \\\hline\hline
			9 & $m-9$  & $m$ = 18 to 30 \\\hline
			10 & $m-10$  & $m$ = 20 to 33 \\\hline\hline
			11 & $m-11$  & $m$ = 22 to 37 \\\hline
			12 & $m-12$  & $m$ = 24 to 40 \\\hline
			13 & $m-13$  & $m$ = 26 to 43 \\\hline\hline
			14 & $m-14$ & $m$ = 28 to 47 \\\hline
			15 & $m-15$  & $m$ = 30 to 50 \\\hline\hline
			16 & $m-16$  & $m$ = 32 to 54 \\\hline
			17 & $m-17$  & $m$ = 34 to 57 \\\hline
			18 & $m-18$  & $m$ = 36 to 60 \\\hline\hline
			19 & $m-19$  & $m$ = 38 to 64 \\\hline
			20 & $m-20$  & $m$ = 40 to 67 \\\hline
			21 & $m-21$  & $m$ = 42 to 70 \\\hline\hline
			22 & $m-22$  & $m$ = 44 to 74 \\\hline
		\end{tabular}
	\end{center}
\end{table} 
%%%%%%%%%%%%%%%%%

\vspace{.3cm}
\noindent
\textbf{Conclusion} It is clear that unlike DM labeling or $\sum$-labeling, $k$-DM labeling covers more families of graphs and has a lot of scope for further research. 

\vspace{.2cm}
\noindent
\textbf{Declaration of competing interest} 
There is no competing interest.

\vspace{.2cm}
\noindent
\textbf{Acknowledgements}\quad 
We express our sincere thanks to the Central University of Kerala, Kasaragod, Kerala for providing facilities to do this research work.

\newpage
%\vspace{.5cm}
\begin{center}
\textbf{ANNEXURE}
\end{center}

\textbf{Calculation of $m$ for a given $m_1$ = $|S_1|$.} 

\begin{enumerate}
\item [\rm 1.] $m_1$ = 4.

When $m_1$ = 4, $m_2$ = $m-4 \geq m_1 = 4$, $m$ = $m_1+m_2 \geq 8$, 

$m+(m-1)+(m-2)+(m-3)$ = $4m-6$ and $\sum^{m-4}_{i=1}{i}$ = $\frac{(m-4)(m-3)}{2}$. 

Let us see how far $\sum^{m_1}_{i=1}{m-i+1} \geq \sum^{m_2}_{i=1}{i}$ for $m_1$ = 4 and $m \geq 8$. See Table 2.
%%%%%%Table-2

\item [\rm 2.] $m_1$ = 5.

When $m_1$ = 5, $m_2$ = $m-5 \geq m_1 = 5$, $m$ = $m_1+m_2 \geq 10$, 

$m+(m-1)+(m-2)+(m-3)+(m-4)$ = $5m-10$ and $\sum^{m-5}_{i=1}{i}$ = $\frac{(m-5)(m-4)}{2}$. 

See how far $\sum^{m_1}_{i=1}{m-i+1} \geq \sum^{m_2}_{i=1}{i}$ for $m_1$ = 5, $m \geq 10$.  See Table 3.
%%%%%%Table-3

\item [\rm 3.] $m_1$ = 6.

When $m_1$ = 6, $m_2$ = $m-6 \geq m_1 = 6$, $m$ = $m_1+m_2 \geq 12$, 

$\sum^{m_1}_{i=1}{m-i+1}$ = $\sum^{6}_{i=1}{m-i+1}$ = $6m-15$ and $\sum^{m_2}_{i=1}{i}$ = $\sum^{m-6}_{i=1}{i}$ = $\frac{(m-6)(m-5)}{2}$. 

See how far $\sum^{m_1}_{i=1}{m-i+1} \geq \sum^{m_2}_{i=1}{i}$ for $m_1$ = 6, $m \geq 12$.  See Table 4.
%%%%%%Table-4

\item [\rm 4.] $m_1$ = 7.

When $m_1$ = 7, $m_2$ = $m-7 \geq m_1 = 7$, $m$ = $m_1+m_2 \geq 14$, 

$\sum^{m_1}_{i=1}{m-i+1}$ = $\sum^{7}_{i=1}{m-i+1}$ = $7m-21$ and $\sum^{m_2}_{i=1}{i}$ = $\sum^{m-7}_{i=1}{i}$ = $\frac{(m-7)(m-6)}{2}$. 

See how far $\sum^{m_1}_{i=1}{m-i+1} \geq \sum^{m_2}_{i=1}{i}$ for $m_1$ = 7,  $m \geq 14$. See Table 5.
%%%%%%Table-5

\item [\rm 5.] $m_1$ = 8.

When $m_1$ = 8, $m_2$ = $m-8 \geq m_1 = 8$, $m$ = $m_1+m_2 \geq 16$, 

$\sum^{m_1}_{i=1}{m-i+1}$ = $\sum^{8}_{i=1}{m-i+1}$ = $8m-28$ and 

$\sum^{m_2}_{i=1}{i}$ = $\sum^{m-8}_{i=1}{i}$ = $\frac{(m-8)(m-7)}{2}$. 

See how far $\sum^{m_1}_{i=1}{m-i+1} \geq \sum^{m_2}_{i=1}{i}$ for $m_1$ = 8,  $m \geq 16$.  See Table 6.
%%%%%%Table-6

\item [\rm 6.] $m_1$ = 9.

When $m_1$ = 9, $m_2$ = $m-9 \geq m_1 = 9$, $m$ = $m_1+m_2 \geq 18$, 

$\sum^{m_1}_{i=1}{m-i+1}$ = $\sum^{9}_{i=1}{m-i+1}$ = $9m-36$ and 

$\sum^{m_2}_{i=1}{i}$ = $\sum^{m-9}_{i=1}{i}$ = $\frac{(m-9)(m-8)}{2}$. 

See how far $\sum^{m_1}_{i=1}{m-i+1} \geq \sum^{m_2}_{i=1}{i}$ for $m_1$ = 9, $m \geq 18$. See Table 7.
%%%%%%Table-7

\item [\rm 7.] $m_1$ = 10.

When $m_1$ = 10, $m_2$ = $m-10 \geq m_1 = 10$, $m$ = $m_1+m_2 \geq 20$, 

$\sum^{m_1}_{i=1}{m-i+1}$ = $\sum^{10}_{i=1}{m-i+1}$ = $10m-45$ and 

$\sum^{m_2}_{i=1}{i}$ = $\sum^{m-10}_{i=1}{i}$ = $\frac{(m-10)(m-9)}{2}$. 

Let us see how far $\sum^{m_1}_{i=1}{m-i+1} \geq \sum^{m_2}_{i=1}{i}$ for $m_1$ = 10 and $m \geq 20$. See Table 8.
%%%%%%Table-8

\item [\rm 8.] $m_1$ = 11.

When $m_1$ = 11, we have to compare $11m-55$ and $\frac{(m-11)(m-10)}{2}$ for $m \geq 22$, $m\in\mathbb{N}$.  See Table 9.
%%%%%%Table-9

\item [\rm 9.] $m_1$ = 12.

When $m_1$ = 12, we have to compare $12m-66$ and $\frac{(m-12)(m-11)}{2}$ for $m \geq 24$, $m\in\mathbb{N}$.  See Table 10.
%%%%%%Table-10

\item [\rm 10.] $m_1$ = 13.

When $m_1$ = 13, we have to compare $13m-78$ and $\frac{(m-13)(m-12)}{2}$ for $m \geq 26$, $m\in\mathbb{N}$.  See Table 11.

%%%%%%Table-11

\item [\rm 11.] $m_1$ = 14.

When $m_1$ = 14, we have to compare $14m-91$ and $\frac{(m-14)(m-13)}{2}$ for $m \geq 28$, $m\in\mathbb{N}$.  See Table 12.
%%%%%%Table-12

\item [\rm 12.] $m_1$ = 15.

When $m_1$ = 15, we have to compare $15m-105$ and $\frac{(m-15)(m-14)}{2}$ for $m \geq 30$, $m\in\mathbb{N}$. See Table 13.
%%%%%%Table-13

\item [\rm 13.] $m_1$ = 16.

When $m_1$ = 16, we have to compare $16m-120$ and $\frac{(m-16)(m-15)}{2}$ for $m \geq 32$, $m\in\mathbb{N}$. See Table 14.
%%%%%%Table-14

\item [\rm 14.] $m_1$ = 17.

When $m_1$ = 17, we have to compare $17m-136$ and $\frac{(m-17)(m-16)}{2}$ for $m \geq 34$, $m\in\mathbb{N}$. See Table 15.
%%%%%%Table-15

\item [\rm 15.] $m_1$ = 18.

When $m_1$ = 18, we have to compare $18m-153$ and $\frac{(m-18)(m-17)}{2}$ for $m \geq 36$, $m\in\mathbb{N}$. See Table 16.
%%%%%%Table-16

\item [\rm 16.] $m_1$ = 19.

When $m_1$ = 19, we have to compare $19m-171$ and $\frac{(m-19)(m-18)}{2}$ for $m \geq 38$, $m\in\mathbb{N}$. See Table 17.
%%%%%%Table-17

%%%%%%%%%%%%%%%%%
\item [\rm 17.] $m_1$ = 20.

When $m_1$ = 20, we have to compare $20m-190$ and $\frac{(m-20)(m-19)}{2}$ for $m \geq 40$, $m\in\mathbb{N}$. See Table 18.
%%%%%%Table-18

\item [\rm 18.] $m_1$ = 21.

When $m_1$ = 21, we have to compare $21m-210$ and $\frac{(m-21)(m-20)}{2}$ for $m \geq 42$, $m\in\mathbb{N}$. See Table 19.
%%%%%%%%%%%%%Table-19

\item [\rm 19.] $m_1$ = 22.

When $m_1$ = 22, we have to compare $22m-231$ and $\frac{(m-22)(m-21)}{2}$ for $m \geq 44$, $m\in\mathbb{N}$. See Table 20.
%%%%%%Table-20

%%%%%%Table-2
\begin{table} \label{t2} 
	\caption{\small{Comparing $4m-6$ and $\frac{(m-4)(m-3)}{2}$ for $m \geq 8$, $m\in\mathbb{N}$.}}
	\begin{center}
		%\scalebox{0.5}{
		\begin{tabular}{|c|c|c|c|} \hline \hline
			\hspace{.2cm} $m$ \hspace{.2cm}
			& $4m-6$  & \hspace{.3cm} $\frac{(m-4)(m-3)}{2}$ \hspace{.3cm} & \hspace{.1cm} $4m-6 \geq\frac{(m-4)(m-3)}{2}$ is True or False \hspace{.1cm}  \\\hline \hline 
			8 & 26  & 10 &  T  \\\hline 
			9 & 30  & 15 &  T  \\\hline 
			10 &  34 & 21 &  T  \\\hline 
			11 & 38  & 28 &  T  \\\hline 
			12 & 42  & 36 &  T  \\\hline 
			13 & 46  & 45 &  T  \\\hline 
			14 & 50  & 55 &  F  \\\hline 
		\end{tabular}
	\end{center}
\end{table} 
%%%%%%Table-3
\begin{table} \label{t3} 
	\caption{\small{Comparing $5m-10$ and $\frac{(m-5)(m-4)}{2}$ for $m \geq 10$, $m\in\mathbb{N}$.}}
	\begin{center}
		%\scalebox{0.5}{
		\begin{tabular}{|c|c|c|c|} \hline \hline
			\hspace{.2cm} $m$ \hspace{.2cm}
			& $5m-10$  & \hspace{.3cm} $\frac{(m-5)(m-4)}{2}$ \hspace{.3cm} & \hspace{.1cm} $5m-10 \geq\frac{(m-5)(m-4)}{2}$ is True or False \hspace{.1cm}  \\\hline \hline 
			10 & 40  & 15 &  T  \\\hline 
			\vdots &\vdots &\vdots  &\vdots  \\\hline 
			15  & 65  & 55 &  T  \\\hline 
			16  & 70  & 66 &  T  \\\hline 
			17  & 75  & 78 &  F  \\\hline 
		\end{tabular}
	\end{center}
\end{table} 
%%%%%%Table-4
\begin{table} \label{t4} 
\caption{\small{Comparing $6m-15$ and $\frac{(m-6)(m-5)}{2}$ for $m \geq 12$, $m\in\mathbb{N}$.}}
\begin{center}
%\scalebox{0.5}{
\begin{tabular}{|c|c|c|c|} \hline \hline
\hspace{.2cm} $m$ \hspace{.2cm}
& $6m-15$  & \hspace{.3cm} $\frac{(m-6)(m-5)}{2}$ \hspace{.3cm} & \hspace{.1cm} $6m-15 \geq\frac{(m-6)(m-5)}{2}$ is True or False \hspace{.1cm}  \\\hline \hline 
12  & 57  & 21 &  T  \\\hline 
\vdots &\vdots &\vdots &\vdots  \\\hline 
19  &  99 & 91 &  T  \\\hline 
20  & 105  & 105 &  T  \\\hline 
21  & 111  & 120 &  F  \\\hline 	
\end{tabular}
\end{center}
\end{table} 
%%%%%%Table-5
\begin{table} \label{t5} 
\caption{\small{Comparing $7m-21$ and $\frac{(m-7)(m-6)}{2}$ for $m \geq 14$, $m\in\mathbb{N}$.}}
\begin{center}
%\scalebox{0.5}{
\begin{tabular}{|c|c|c|c|} \hline \hline
\hspace{.2cm} $m$ \hspace{.2cm}
& $7m-21$  & \hspace{.3cm} $\frac{(m-7)(m-6)}{2}$ \hspace{.3cm} & \hspace{.1cm} $7m-21 \geq\frac{(m-7)(m-6)}{2}$ is True or False \hspace{.1cm}  \\\hline \hline 
14  & 77 & 28 &  T  \\\hline 
\vdots &\vdots &\vdots &\vdots  \\\hline 
22  & 133  & 120 &  T  \\\hline 
23  & 140  & 136 &  T  \\\hline 
24  & 147  & 153 &  F  \\\hline  	
\end{tabular}
\end{center}
\end{table} 
%%%%%%Table-6
\begin{table} \label{t6} 
\caption{\small{Comparing $8m-28$ and $\frac{(m-8)(m-7)}{2}$ for $m \geq 16$, $m\in\mathbb{N}$.}}
\begin{center}
%\scalebox{0.5}{
\begin{tabular}{|c|c|c|c|} \hline \hline
\hspace{.2cm} $m$ \hspace{.2cm}
& $8m-28$  & \hspace{.3cm} $\frac{(m-8)(m-7)}{2}$ \hspace{.3cm} & \hspace{.1cm} $8m-28 \geq\frac{(m-8)(m-7)}{2}$ is True or False \hspace{.1cm}  \\\hline \hline 
16  & 100 & 36 &  T  \\\hline 
\vdots &\vdots  &\vdots & \vdots  \\\hline 
 25  & 172 & 153 &  T  \\\hline 
26  & 180 & 171 &  T  \\\hline 
27  & 188 & 190 &  F  \\\hline 
\end{tabular}
\end{center}
\end{table} 
%%%%%%Table-7
\begin{table} \label{t7} 
\caption{\small{Comparing $9m-36$ and $\frac{(m-9)(m-8)}{2}$ for $m \geq 18$, $m\in\mathbb{N}$.}}
\begin{center}
%\scalebox{0.5}{
\begin{tabular}{|c|c|c|c|} \hline \hline
\hspace{.2cm} $m$ \hspace{.2cm}
& $9m-36$  & \hspace{.3cm} $\frac{(m-9)(m-8)}{2}$ \hspace{.3cm} & \hspace{.1cm} $9m-36 \geq\frac{(m-9)(m-8)}{2}$ is True or False \hspace{.1cm}  \\\hline \hline 
18  & 126 & 45 &  T  \\\hline 
\vdots & \vdots  & \vdots  &  \vdots   \\\hline 
29  & 225 & 210 &  T  \\\hline 
30  & 234 & 231 &  T  \\\hline 
31  & 243 & 253 &  F  \\\hline 
\end{tabular}
\end{center}
\end{table} 
%%%%%%Table-8
\begin{table} \label{t8}
\caption{\small{Comparing $10m-45$ and $\frac{(m-10)(m-9)}{2}$ for $m \geq 20$, $m\in\mathbb{N}$.}}
\begin{center}
%\scalebox{0.5}{
\begin{tabular}{|c|c|c|c|} \hline \hline
 $m$ & $10m-45$  & $\frac{(m-10)(m-9)}{2}$  &  $10m-45 \geq\frac{(m-10)(m-9)}{2}$ is True or False   \\\hline \hline 
20  & 155 & 55 &  T  \\\hline 
\vdots & \vdots & \vdots & \vdots \\\hline 
32  & 275 & 253 &  T  \\\hline 
33  & 285 & 276 &  T  \\\hline 
34  & 295 & 300 &  F  \\\hline 
\end{tabular}
\end{center}
\end{table} 
%%%%%%Table-9
\begin{table} \label{t9} 
\caption{\small{Comparing $11m-55$ and $\frac{(m-11)(m-10)}{2}$ for $m \geq 22$, $m\in\mathbb{N}$.}}
\begin{center}
%\scalebox{0.5}{
\begin{tabular}{|c|c|c|c|} \hline \hline
 $m$ & $11m-55$  &  $\frac{(m-11)(m-10)}{2}$  &  $11m-55 \geq\frac{(m-11)(m-10)}{2}$ is True or Fales   \\\hline \hline 
22  & 187 & 66 &  T  \\\hline 
\vdots & \vdots & \vdots & \vdots \\\hline 
36  & 341 & 325 &  T  \\\hline 
37  & 352 & 351 &  T  \\\hline 
38  & 363 & 378 &  F  \\\hline 
\end{tabular}
\end{center}
\end{table} 
%%%%%%Table-10
\begin{table} \label{t10} 
\caption{\small{Comparing $12m-66$ and $\frac{(m-12)(m-11)}{2}$ for $m \geq 24$, $m\in\mathbb{N}$.}}
\begin{center}
%\scalebox{0.5}{
\begin{tabular}{|c|c|c|c|} \hline \hline
 $m$ & $12m-66$  &  $\frac{(m-12)(m-11)}{2}$  &  $12m-66 \geq\frac{(m-12)(m-11)}{2}$ is True or False   \\\hline \hline 
24  & 222 & 78 &  T  \\\hline 
\vdots & \vdots & \vdots & \vdots \\\hline 
40  & 414 & 406 &  T  \\\hline 
41  & 426 & 435 &  F  \\\hline 
\end{tabular}
\end{center}
\end{table} 
%%%%%%Table-11
\begin{table} \label{t11} 
\caption{\small{Comparing $13m-78$ and $\frac{(m-13)(m-12)}{2}$ for $m \geq 26$, $m\in\mathbb{N}$.}}
\begin{center}
%\scalebox{0.5}{
\begin{tabular}{|c|c|c|c|} \hline \hline
$m$ & $13m-78$  &  $\frac{(m-13)(m-12)}{2}$  &  $13m-78 \geq\frac{(m-13)(m-12)}{2}$ is True or Fales  \\\hline \hline 
26  & 260 & 91 &  T  \\\hline 
\vdots & \vdots & \vdots & \vdots \\\hline 
43  & 481 & 465 &  T  \\\hline 
44  & 494 & 496 &  F  \\\hline 
\end{tabular}
\end{center}
\end{table} 
%%%%%%Table-12
\begin{table} \label{t12} 
\caption{\small{Comparing $14m-91$ and $\frac{(m-14)(m-13)}{2}$ for $m \geq 28$, $m\in\mathbb{N}$.}}
\begin{center}
%\scalebox{0.5}{
\begin{tabular}{|c|c|c|c|} \hline \hline
 $m$ & $14m-91$  &  $\frac{(m-14)(m-13)}{2}$ &  $14m-91 \geq\frac{(m-14)(m-13)}{2}$ is True or False   \\\hline \hline 
28  & 301 & 105 &  T  \\\hline 
\vdots & \vdots & \vdots & \vdots \\\hline 
47  & 567 & 561 &  T  \\\hline 
48  & 581 & 595 &  F  \\\hline 
\end{tabular}
\end{center}
\end{table} 
%%%%%%Table-13
\begin{table} \label{t13} 
\caption{\small{Comparing $15m-105$ and $\frac{(m-15)(m-14)}{2}$ for $m \geq 30$, $m\in\mathbb{N}$.}}
\begin{center}
%\scalebox{0.5}{
\begin{tabular}{|c|c|c|c|} \hline \hline
 $m$ & $15m-105$  & $\frac{(m-15)(m-14)}{2}$ &  $15m-105 \geq\frac{(m-15)(m-14)}{2}$ is True or False  \\\hline \hline 
30  & 345 & 120 &  T  \\\hline 
\vdots & \vdots & \vdots & \vdots \\\hline 
50  & 645 & 630 &  T  \\\hline 
51  & 660 & 666 &  F  \\\hline 
\end{tabular}
\end{center}
\end{table} 
%%%%%%Table-14
\begin{table} \label{t14} 
\caption{\small{Comparing $16m-120$ and $\frac{(m-16)(m-15)}{2}$ for $m \geq 32$, $m\in\mathbb{N}$.}}
\begin{center}
%\scalebox{0.5}{
\begin{tabular}{|c|c|c|c|} \hline \hline
 $m$ & $16m-120$  &  $\frac{(m-16)(m-15)}{2}$  & $16m-120 \geq\frac{(m-16)(m-15)}{2}$ is True or False   \\\hline \hline 
32  & 392 & 136 &  T  \\\hline 
\vdots & \vdots & \vdots & \vdots \\\hline 
53  & 728 & 703 &  T  \\\hline 
54  & 744 & 741 &  T  \\\hline 
55  & 760 & 780 &  F  \\\hline 
\end{tabular}
\end{center}
\end{table} 
%%%%%%Table-15
\begin{table} \label{t15} 
\caption{\small{Comparing $17m-136$ and $\frac{(m-17)(m-16)}{2}$ for $m \geq 34$, $m\in\mathbb{N}$.}}
\begin{center}
%\scalebox{0.5}{
\begin{tabular}{|c|c|c|c|} \hline \hline
 $m$ 
& $17m-136$  &  $\frac{(m-17)(m-16)}{2}$  &  $17m-136 \geq\frac{(m-17)(m-16)}{2}$ is True or False   \\\hline \hline 
34  & 442 & 153 &  T  \\\hline 
\vdots & \vdots & \vdots & \vdots \\\hline 
57  & 833 & 820 &  T  \\\hline 
58  & 850 & 861 &  F  \\\hline 
\end{tabular}
\end{center}
\end{table} 
%%%%%%Table-16
\begin{table} \label{t16} 
\caption{\small{Comparing $18m-153$ and $\frac{(m-18)(m-17)}{2}$ for $m \geq 36$, $m\in\mathbb{N}$.}}
\begin{center}
%\scalebox{0.5}{
\begin{tabular}{|c|c|c|c|} \hline \hline
 $m$ & $18m-153$  &  $\frac{(m-18)(m-17)}{2}$  &  $18m-153 \geq\frac{(m-18)(m-17)}{2}$ is True or False   \\\hline \hline 
36  & 495 & 171 &  T  \\\hline 
\vdots & \vdots & \vdots & \vdots \\\hline 
60  & 927 & 903 &  T  \\\hline 
61  & 945 & 946 &  F  \\\hline 
\end{tabular}
\end{center}
\end{table} 
%%%%%%Table-17
\begin{table} \label{t17} 
\caption{\small{Comparing $19m-171$ and $\frac{(m-19)(m-18)}{2}$ for $m \geq 38$, $m\in\mathbb{N}$.}}
\begin{center}
%\scalebox{0.5}{
\begin{tabular}{|c|c|c|c|} \hline \hline
 $m$ & $19m-171$  &  $\frac{(m-19)(m-18)}{2}$  &  $19m-171 \geq\frac{(m-19)(m-18)}{2}$ is True or False   \\\hline \hline 
38  &  & 190 &  T  \\\hline 
\vdots & \vdots & \vdots & \vdots \\\hline 
63  & 1026 & 990 &  T  \\\hline 
64  & 1045 & 1035 &  T  \\\hline 
65  & 1064 & 1081 &  F  \\\hline 
\end{tabular}
\end{center}
\end{table} 
%%%%%%Table-18
\begin{table} \label{t18} 
\caption{\small{Comparing $20m-190$ and $\frac{(m-20)(m-19)}{2}$ for $m \geq 40$, $m\in\mathbb{N}$.}}
\begin{center}
%\scalebox{0.5}{
\begin{tabular}{|c|c|c|c|} \hline \hline
 $m$ & $20m-190$  &  $\frac{(m-20)(m-19)}{2}$  &  $20m-190 \geq\frac{(m-20)(m-19)}{2}$ is True or False   \\\hline \hline 
40  & 610 & 210 &  T  \\\hline 
\vdots & \vdots & \vdots & \vdots \\\hline 
67  & 1150 & 1128 &  T  \\\hline 
68  & 1170 & 1176 &  F  \\\hline 
\end{tabular}
\end{center}
\end{table} 
%%%%%%Table-19
\begin{table} \label{t19} 
\caption{\small{Comparing $21m-210$ and $\frac{(m-21)(m-20)}{2}$ for $m \geq 42$, $m\in\mathbb{N}$.}}
\begin{center}
%\scalebox{0.5}{
\begin{tabular}{|c|c|c|c|} \hline \hline
 $m$ & $21m-210$  &  $\frac{(m-21)(m-20)}{2}$  &  $21m-210 \geq\frac{(m-21)(m-20)}{2}$ is True or False   \\\hline \hline 
42  & 672 & 231 &  T  \\\hline 
\vdots & \vdots & \vdots & \vdots \\\hline 
70  & 1260 & 1225 &  T  \\\hline 
71  & 1281 & 1275 &  F  \\\hline 
\end{tabular}
\end{center}
\end{table} 
%%%%%%Table-20
\begin{table} \label{t20} 
\caption{\small{Comparing $22m-231$ and $\frac{(m-22)(m-21)}{2}$ for $m \geq 44$, $m\in\mathbb{N}$.}}
\begin{center}
%\scalebox{0.5}{
\begin{tabular}{|c|c|c|c|} \hline \hline
 $m$ & $22m-231$  &  $\frac{(m-22)(m-21)}{2}$  &  $22m-231 \geq\frac{(m-22)(m-21)}{2}$ is True or False   \\\hline \hline 
44  & 737 & 253 &  T  \\\hline 
\vdots & \vdots & \vdots & \vdots \\\hline 
74  & 1397 & 1378 &  T  \\\hline 
75  & 1419 & 1431 &  F  \\\hline 
\end{tabular}
\end{center}
\end{table} 
%%%%%%%%%%%%%%%%%
\end{enumerate}

\end{document}